\documentclass[journal]{IEEEtran}

\usepackage[pdftex]{graphicx}
\DeclareGraphicsExtensions{.pdf,.jpeg,.png}
\usepackage[cmex10]{amsmath}
\usepackage{array,color}
  \usepackage[caption=false,font=footnotesize]{subfig}

\usepackage{amsmath} % Needs to be defined before amsthm package
\usepackage{amssymb,url}
\usepackage{amsmath,amsthm, amssymb,url,amsfonts,mathrsfs}

\renewcommand{\v}[1]{\ensuremath{\boldsymbol{\mathbf{#1}}}}

\long\def\old#1{}

\newcommand{\Rset}{\mathbb{R}}
\newcommand{\bdashlist}{\begin{list}{$\bullet$}{} }
\newcounter{l3}
\newcounter{l1}
\newcommand{\bromalist}{\begin{list}{(\roman{l3})}{\usecounter{l3}}}
\newcommand{\barablist}{\begin{list}{\arabic{l1}.}{\usecounter{l1}}}
\newcommand{\Scal}{{\cal S}}
\newcommand{\Nset}{\mathbb{N}}
\newcommand{\Acal}{{\cal A}}
\newcommand{\Pcal}{{\cal P}}
\newcommand{\Zcal}{{\cal Z}}
\newcommand{\Ucal}{{\cal U}}

\newcommand{\Rp}{\mathbb{R}_{+}}

\newcommand{\E}{\mathbb{E}}

\newcommand{\bx}{\mathbf{x}}

\theoremstyle{definition}
\newtheorem{theorem}{Theorem}
\newtheorem{lemma}{Lemma}
\newtheorem{assumption}{Assumption}
\newtheorem{example}{Example}
\newtheorem{definition}{Definition}
\newtheorem{remark}{Remark}
\newtheorem{corollary}{Corollary}

\long\def\blue#1{{\color{black}#1}}
\long\def\yunjian#1{#1}

\begin{document}
%\begin{titlepage}

\title{Deadline Differentiated Pricing of Deferrable Electric Loads}

\author{Eilyan Bitar$^\dagger$ \ \ and \ \ Yunjian Xu$^{\dagger \dagger}$ %<-this % stops a space
\thanks{This work was supported in part by  NSF ECCS-1351621, NSF CNS-1239178, NSF IIP-1632124,  PSERC under sub-award S-52, US DoE under the CERTS initiative, and the MIT-SUTD International Design Center (IDC) Grant IDG21400103. Both authors contributed equally to the submitted work.}% <-this % stops a space
\thanks{$\dagger$ E. Bitar is with the School of Electrical and Computer Engineering, and the School of Operations Research and Information Engineering, Cornell University, Ithaca, NY, 14853, USA.  Email:  {\tt\small eyb5@cornell.edu}}
\thanks{$\dagger\dagger\,$Y. Xu is with the Engineering Systems and Design Pillar, Singapore University of Technology and Design, Singapore.   Email:  {\tt\small yunjian\_xu@sutd.edu.sg}}
%\thanks{The authors would like to thank Costas
%Courcoubetis, Stephen Graves, Steven Low, Paul De Martini, Kameshwar Poolla,
%Pravin Varaiya, and Adam Wierman for their helpful discussions}
}

\maketitle

\begin{abstract}
A large fraction of total electricity demand is comprised of end-use devices whose demand for energy is inherently deferrable in time. Of interest is the potential to use this latent flexibility in demand to absorb variability in power supplied from intermittent renewable generation.
A fundamental challenge lies in the design of incentives that induce the desired response in demand. With an eye to electric vehicle charging, we propose a novel forward market for deadline-differentiated electric power service, where consumers consent to deferred service of pre-specified loads in exchange
for a reduced price for energy. The longer a consumer is willing to defer, the lower the price for energy.
The proposed forward contract provides a guarantee on the aggregate quantity of energy to be delivered by a
consumer-specified deadline. Under the earliest-deadline-first (EDF) scheduling policy, which is shown to be optimal for the supplier, we explicitly characterize a non-discriminatory, deadline-differentiated pricing scheme that yields an efficient competitive equilibrium between the supplier and consumers. We further show
that this efficient pricing scheme, in combination with EDF scheduling, is incentive compatible (IC) in that every consumer would like to reveal her true deadline to the supplier, regardless of the actions taken by other consumers.
\end{abstract}
%\thispagestyle{empty}
%\end{titlepage}

%\vspace{-.15in}

\IEEEpeerreviewmaketitle

\vspace{-.15in}
%%%%%%%%%%%%%%%%%%%%%%%%%%%%%%%%%%%%%%%%%%%%%%%%%%%%%%%%%%%%%%%%%%%%%%%%%
\section{Introduction} \label{sec:intro}

\noindent As the electric power industry transitions to a greater
reliance on intermittent and distributed energy resources, there emerges a need for flexible resources that can respond
dynamically to weather impacts on wind and solar photovoltaic
output. These renewable generation sources have limited
controllability and production patterns that are intermittent and
uncertain. Such variability represents one \yunjian{of} the most important
obstacles to the deep integration of renewable generation into the
electricity grid. The current approach to renewable energy
integration is to balance variability with dispatchable generation.
This works at today's modest penetration levels, but it cannot
scale, because of the projected increase in reserve generation
required to balance the variability in renewable supply
\cite{CAISO2010}.
%Recent studies in California \cite{CAISO2010} project that the spring-time maximum up-regulation capacity needed to accommodate 33\% renewable energy penetration will increase from 277 MW to 1,135 MW. Similar increases are projected in down-regulation capacity, and maximum load-following capacity requirements will need to increase from 2,292 MW to 4,423 MW.
If these increases are met with combustion fired generation, they
will both be counterproductive to carbon emissions reductions and
economically untenable.

As wind and solar energy penetration increases, how must the
assimilation of this variable power evolve, so as to minimize these
integration costs, while maximizing the net environmental benefit?
Clearly, strategies which attenuate the increase in conventional
reserve requirements will be an essential means to this end. One
option is to harness the flexibility in demand-side resources.
As such, significant benefits have been identified by the
Federal Energy Regulatory Commission \cite{FERC2009} in
unlocking the value of coordinating demand-side resources to
address the growing need for firm, responsive resources to provide
 balancing services for the bulk
power system.

\subsection{The Current Approach to Demand Response}
\noindent There is an opportunity to transform the current
operational paradigm,
 in which supply is tailored to follow demand, to one in which \emph{demand is capable of reacting to variability in supply} --
 an approach which is generally referred to as demand response (DR)
 \cite{Albadi07}.
The primary challenge is the \emph{reliable extraction of the desired response}
 from participating demand resources on time scales aligned with traditional bulk power balancing services.

The majority of DR programs in place today are limited to peak
shaving and contingency-based applications. The two most common paradigms
for customer recruitment and control are: (1) \emph{direct load
control} where a load aggregator or utility procures the capability of direct load adjustment through a forward market or transaction
%\blue{under which the system operator or utility companies procure the capability of load adjustment through a forward transaction}
(e.g. call options for interruptible load) \cite{Chao87, NN14, Oren92, Varaiya92} and (2) \emph{indirect load
control} where consumers or devices adjust their energy consumption
 in response to dynamic  price signals \cite{Fuller11, H10, LiChen11, LiZhang15}.
While dynamic pricing has the  potential to improve the economic efficiency of electricity markets
\cite{B05,SL07,IX11}, it subjects consumers to the risk of
paying high peak prices, and can have the counterproductive effect of increasing the variability in demand \cite{MCK11, RD12}.
And, of particular relevance to this paper's emphasis on electric vehicle (EV) charging is a recent empirical study \cite{S12} that
indicates dynamic pricing may perform worse than a
flat-rate tariff for
EV charging in terms of generation costs and emissions impacts.
%
%A recent market monitoring report
%expresses the concern of California Independent System Operator (CAISO) on dynamic  pricing: ``While there are many economists that are enthusiastic about DR for all consumers, we are not aware of a reported success of real-time pricing for a big, heterogeneous population area that could serve as a benchmark. Mobilizing retail level demand side flexibility to reduce operating and investment cost in the electricity sector by employing smart grid technologies and market mechanism is still regarded as work in progress.'' \red{(REFERENCE ??)}
%
%
In short, demand response implemented through dynamic pricing may not provide the level of assurance required to avoid the use of conventional
generation to manage the electric power system.

%\begin{quote}
%{\it ``While there are many economists that are enthusiastic about DR for all consumers, we are not aware of a reported success of real-time pricing for a big, heterogeneous population area that could serve as a benchmark. Mobilizing retail level demand side flexibility to reduce operating and investment cost in the electricity sector by employing smart grid technologies and market mechanism is still regarded as work in progress.''}
%\end{quote}

\old{As the use of responsive demand shifts from reliability based
utility run programs to market-based real-time balancing services,
the criteria for performance guarantees become more stringent.}

\subsection{A Deadline Differentiated Energy Service Approach}
In the following paper, we propose a novel market framework  to enable the direct control of deferrable loads, with a particular emphasis on electric vehicle (EV) charging. Broadly, the proposed market centers on the provisioning of deadline differentiated energy services to customers possessing the inherent ability to delay their consumption (up to a point) in time.
%\footnote{Such
%a market would naturally complement the increased proliferation of
%plug-in electric vehicles (EV) in the US transportation fleet.}
From the consumer's perspective, the longer she is willing to delay the receipt of a specified quantity of
energy, the less that customer pays (per-unit) for said energy.
The supplier, on the other hand, implicitly purchases the right to manage the real-time
delivery of power to participating consumers by offering a discount on energy with longer deadlines on delivery. And, the longer the consumer-specified deadlines, the more flexibility the supplier has in meeting the corresponding energy requirements. Put simply, the supplier can extract flexibility from the demand-side through direct load control, all while providing firm guarantees on delivery by consumer-specified deadlines. In this way, the supplier can
align its operational requirements with the vast heterogeneity in
end-use customer needs.

%
%
%
%\subsection{A Deadline Differentiated Service Approach}
%\noindent Flexibility in consumption can be interpreted as a
%\emph{continuum of feasible power profiles} capable of preserving
%the end-use function of a demand resource. A basic question, is how
%to design a market that enables a coordinating entity the ability to
%``extract this flexibility'' for execution of real-time control
%applications -- e.g. balancing variability in renewable supply?
%
%
%
%One possibility resides in the construction of a market for
%\emph{quality-differentiated electric power services}, where the
%price to a consumer for receiving a particular service is a
%monotonic function of the desired quality-of-service (QoS).
%Naturally, a reduction in QoS is accompanied by a reduction in
%price.  For example, in the concrete setting of deferrable loads,
%deadline would be a natural specification of QoS. The longer a
%customer is willing to delay the receipt of a specified quantity of
%energy, the less that customer pays (per-unit) for said energy.\old{Such markets would require a shift in how we think of electricity --
%not as an undifferentiated good, but rather as a set of
%differentiated services from which individual consumers can specify
%a QoS that best meets their electricity needs. Moreover,} As a QoS
%specification maps directly to a set of feasible power profiles
%capable of servicing a load, the aggregator can imbed its extraction
%of flexibility from the demand-side in its delivery of
%differentiated services with a guaranteed QoS to participating
%consumers.

The general concept of electric power service differentiation is not new
\cite{Oren92,Oren93}. Many have studied the problem of centralized control of a collection of loads for
load-following or regulation services -- all while ensuring the
satisfaction of a pre-specified quality-of-service (QoS) to individual resources
\cite{Callaway11, HaoChen14, LangAllerton2011, NM12, Obrien13, Papa2010PESGM,
Roozbehani2011, Sub2012ACC, YC14}.
There has, however, been little
work in the way of designing market mechanisms that endogenously
price the flexibility being offered by the demand side, while
incentivizing consumers to truthfully reveal their preferences  to the operator -- for example,
the ability to delay energy consumption in time. Along these lines, the authors in \cite{NN14} propose a forward market
for duration-differentiated energy services. Several classic \cite{Chao87, Varaiya92} and more recent \cite{Bitar2012HICSS} papers have explored the concept
of reliability-differentiated pricing  of interruptible electric
power service, where consumers take on the risk of supply interruption in exchange for a reduction in the nominal energy price.
\old{Similar in spirit, \cite{YC14} devise a market for \emph{risk-limiting contracts}, where each consumer can explicitly specify the maximal price risk he or she is willing to accept.}
Beyond the
 difficulty in auditing such markets and the apparent issues of moral hazard, the primary drawback of such
approaches stems from the explicit transferal of quantity or price risk to
the demand side. This amounts to requiring that consumers plan their consumption in the face of uncertain supply or prices -- a nontrivial decision task.

%as it requires participating consumers to plan
%their consumption in the face of uncertain supply. From a consumer's perspective, this amounts to solving a nontrivial %

With the aim of alleviating the aforementioned challenges,
we propose and analyze a novel forward market for \emph{deadline differentiated energy
services},\footnote{\blue{The concept of deadline differentiated pricing was first proposed
in a conference paper \cite{BitarLow2012} and subsequently empirically evaluated in \cite{SF14}.
%In this paper, we provide characterization on the efficient competitive equilibrium at which every consumer is truth-telling.
}}  where consumers consent to deferred supply of energy in exchange for a reduced per-unit {energy price}.
Such a market would naturally complement the proliferation of
plug-in electric vehicles (EV) in the US transportation fleet.

\begin{example}[Electric Vehicle Charging]
With this motivation in mind, we now illustrate several basic questions that might arise in the design and operation of markets for deadline differentiated energy services. To frame the discussion, consider a scenario involving the operation of a parking infrastructure equipped with an array of EV charging stations and local photovoltaic (PV) generation.\footnote{Of particular relevance to our development is the growing multitude of companies offering turn-key products that integrate EV charging infrastructure with solar photovoltaic canopies \cite{EVSO14, SolGen11}.} The EVs being charged in such facilities  are naturally modeled as having the ability to delay their receipt of energy in time. Accordingly, upon connecting an EV to a  charging station, the vehicle's owner is presented with
a menu of prices --  each of which stipulates a per-unit
price for energy and a corresponding delivery deadline.
Faced with such choices, how does the EV owner decide upon which
bundle of energy-deadline pairs to purchase given the inherent ability to delay consumption?
On the supply side, the operator of the parking infrastructure is required to meet all energy requests by their corresponding deadlines. How should the operator set the menu of deadline differentiated
prices to induce the desired demand from a population of EV owners?
What if the available energy supply is random, as is the case with PV generation? In order to address such questions, we first require a description of the underlying market, which we informally describe below to orient the reader.  $\hfill \Box$
%An EV owner arrives at the car park and connects her EV to charging station.
%Upon plugging in, the EV owner is presented with
%a menu of prices --  each of which stipulates a per-unit
%price for energy and a corresponding delivery deadline.
%Faced with such choices, how does the EV owner decide upon which
%bundle of energy-deadline pairs to purchase given the inherent ability to delay consumption?
%On the supply side, the load aggregator is required to meet all energy requests by their
%corresponding deadlines. How should the aggregator set the menu of deadline differentiated
%prices to induce the desired demand from a population of EV owners?
%What if the available energy supply is random? In order to address such questions, we first require a description of the underlying market, which we informally describe below to orient the reader.  $\hfill \Box$
\end{example}

\noindent \textbf{Market Operation.} \  The forward market for deadline differentiated energy
service operates according to a three-step process. In step 1, the supplier announces a mechanism. In step 2, the consumers simultaneously report their demand. In step 3, the mechanism is executed; namely, prices are set and the requested demand is delivered to each consumer. Time is assumed to be
discrete with periods indexed by $k=0,1,2,\dots,N$.

\

\noindent \textbf{Step 1} (Mechanism Design). $\;$ \yunjian{At the beginning of} period $k=0$,
the supplier announces a \emph{market mechanism} $(\pi,\v{\kappa})$ consisting of both a
\emph{scheduling policy} $\pi$  (cf. its formal definition in Section \ref{sec:supply_model}),
and  \emph{pricing scheme} $\v{\kappa}=(\kappa_1,\ldots,\kappa_N)$ that maps the aggregate demand bundle $\bx$ \ (cf. its definition in Step 2) into a
menu of \emph{deadline-differentiated
prices},
%\begin{align*}
%\mathcal{M} = \left\{  (p_k, d_k) \in \Rset_+^2 \ | \  0 \leq d_0 < d_1 < \cdots < d_M  \right\},
%\end{align*}
\begin{equation}\label{eq:ka}
p_k=\kappa_k(\bx), \qquad k=1,\ldots,N.
\end{equation}
The price menu stipulates a per-unit price $p_k$ (\$/kWh) for energy guaranteed
delivery by period $k$. At the heart of the mechanism design is the restriction that prices are nonincreasing in the deadline.
%$$
%p_1 \ \ge \ p_2  \ \ge \ldots \ge \ p_N.
%$$
Namely, the longer a customer is willing to defer her consumption, the less she is required to pay.
We will use $\mathcal{P}  = \{ \v{p} \in \Rset_+^N \ | \ p_1  \ge  p_2   \ge \ldots \ge  p_N\} $ to denote the set of
feasible deadline-differentiated price bundles.

\

\noindent \textbf{Step 2} (Consumer Reporting). $\;$  Each consumer then
reports a bundle of deadline-differentiated energy quantities $\mathbf{a} = (a_1,\dots, a_N)^{\top} \in
\Rp^{N}$. Here, the quantity
$a_k$ (kWh) denotes the amount of energy guaranteed delivery by the deadline $k$. It follows that said consumer will receive at least $\sum_{t=1}^k a_t$ amount of energy by deadline $k$.
%The consumer is required to pay
%$ \v{p}^{\top} \v{a}$.
The aggregate demand bundle is the sum of all individual consumer bundles, which we denote by  $\mathbf{x} =  (x_1, \dots,
x_N)^{\top} \in \Rset_+^{N}$.
%It follows that $x_k$ denotes the aggregate
%quantity the supplier is required to deliver by deadline $k$.

%
%$\bx$, according to the pricing scheme $\v{\kappa}$ announced prior to period $k=0$. A pricing scheme $\v{\kappa}$ is therefore a mapping from $\Rset_+^{N}$ to $\mathcal{P}$.

\

\noindent \textbf{Step 3} (Pricing and Energy Delivery). $\;$  Given an aggregate demand bundle $\v{x}$, the
deadline differentiated prices are set according to $\v{p} = \v{\kappa}(\v{x})$. Pricing is non-discriminatory, in that all customers are charged according to the same menu of prices. Thus, a customer requesting a bundle $\v{a}$ pays $\v{p}^\top \v{a}$.
The supplier
must also deliver the aggregate demand bundle $\mathbf{x}$ according to the previously announced scheduling policy
$\pi$.  Essentially,  a scheduling policy is said to be feasible if it delivers each
 consumer's requested energy bundle by its corresponding deadline. The
supplier is assumed to have two sources of electricity from
which \yunjian{she} can service demand: \ \emph{intermittent} and \emph{firm}.

\begin{itemize} \setlength{\itemsep}{.1in}
\item \emph{Intermittent supply:} an intermittent supply
modeled as a discrete time random process \ $\mathbf{s} = (s_0, s_1,
\dots, s_{N-1})$ \ with known joint probability distribution.
Here, $s_k \in \Scal$ (kWh) denotes the energy produced during period
$k$ and $\Scal \subset \Rset_+$  the feasible supply interval.
The intermittent supply is assumed to have zero marginal cost.

\item \emph{Firm supply:} a firm supply with a fixed price of $c_0 > 0$. The price $c_0$ can be interpreted as the nominal flat rate for electricity set by the local utility.
\end{itemize}

While stylized in nature, our models of supply and demand are meant to reflect the essential features of an emerging
grid infrastructure that enables the direct coupling of deferrable electric loads with variable renewable supply. 
%Of particular relevance to our development is the growing multitude of companies offering turn-key products that integrate  EV charging infrastructure with solar photovoltaic canopies \cite{EVSO14, SolGen11}.

\subsection{Summary of Main Results}

The primary contribution of our paper is the design of a market mechanism for deadline differentiated energy services, which implements truth-telling at a dominant strategy equilibrium across the population of consumers,
while maximizing social welfare at a competitive equilibrium between a price-taking supplier\footnote{A price-taking supplier cannot control
 market prices, because of government regulation or perfect competition. Instead, it seeks to maximize its profit by scheduling the purchase of its firm supply.} and
 a continuum of infinitesimally small consumers.
We provide here a roadmap of the paper together
with a summary of our main results.

\begin{itemize} \setlength{\itemsep}{.06in}
\item We provide a stylized, yet descriptive,  model of both
deferrable electricity demand
and a supplier with both firm and intermittent  supply
 in Sections \ref{sec:demand_model} and \ref{sec:supply_model}, respectively.

\item We formulate the supplier's scheduling problem as a constrained stochastic optimal control problem
and prove average-cost  optimality of the earliest-deadline-first (EDF) scheduling policy.
% in terms of minimizing the expected cost of firm supply over all feasible scheduling policies.
As a corollary, this result enables
the explicit characterization of the supplier's marginal cost curve in Theorem
\ref{thm:opt_sup}.
%which specifies the bundle of deadline-differentiated
%energy quantities  $\v{x}$ that said supplier would like to provide at a given price
%bundle $\v{p}$.

\item  It is reasonable to expect that the supplier
cannot observe the true deadline of each individual consumer.  The
presence of asymmetric information may lead to significant welfare
loss, if consumers misreport their true deadlines.
%
%
%\red{We note that, under deadline-differentiated pricing,
%a consumer may have incentive to report a false deadline (past her true deadline), if the prices associated with {later} deadlines are low enough.
%In particular, a rational (expected-payoff maximizing) consumer
%will report a false later deadline, if the corresponding saving in energy cost exceeds her expected loss of utility (resulting from the possible
%shortfall in the energy received by her true deadline).}
%
%
Somewhat surprisingly, we show under mild assumptions on consumers' utility functions that a mechanism consisting of an EDF scheduling policy,
together with a uniform marginal cost pricing scheme implements consumers'
\emph{truth-telling} behavior in dominant strategies.

\item   We show in Section \ref{sec:imp} that marginal cost pricing, in combination with EDF scheduling,
induces an \emph{efficient competitive equilibrium}  that simultaneously maximizes the
(price-taking) supplier's profit and the social welfare.
%
%
%\red{In other
%words, the pricing scheme characterized in Theorem
%\ref{thm:opt_sup}, which has been shown to be incentive compatible in
%Section \ref{sec:IC}, results in an efficient market equilibrium, at
%which the social welfare (the sum of aggregate consumer surplus and
%supplier profit) is maximized.
%%
%We also discuss briefly the
%implementation of the proposed deadline differentiated pricing
%mechanism, through both a mechanism design and a (competitive)
%market equilibrium approach.}

\end{itemize}

A preliminary version of this paper appeared in \cite{BX13}, where we  proved the existence of a truth-telling Nash equilibrium.
In the present paper, we establish a stronger
 incentive compatibility result. Namely, 
it is a dominant strategy for each consumer to be truth-telling, regardless of the actions taken by other consumers.
All formal proofs of our stated results can be found in the Appendix of this paper.

\subsection{Related Work}

There have  recently emerged several papers concerned with the design of incentives for deferrable electric loads. The authors of \cite{Kefayati2011} propose an idea similar to this paper, where
consumers are offered a discounted electricity price in exchange for the
delay of their energy consumption. However, the focus of \cite{Kefayati2011} is not on pricing, but rather
on the problem of optimal scheduling faced by the operator, who is assumed to
have full information about consumers -- for example, knowledge of their true deadlines.
Closer to the present paper, \cite{RS12} proposes a greedy online mechanism for electric vehicle charging, which is shown to be incentive compatible
and achieve a bounded (worst-case) competitive ratio.  These theoretical results require, however, a VCG-type (discriminatory-price) payment scheme, as opposed to the uniform-price scheme proposed in this paper.
More strongly, they impose an additional assumption that consumers cannot report false deadlines exceeding their true underlying deadlines. Such assumption substantially simplifies the problem of designing an incentive compatible pricing scheme when consumers are permitted to report arbitrary deadlines -- the setting considered in the present paper.

\section{Model of Demand} \label{sec:demand_model}

\def\figwid{1.68in}
\setlength\fboxsep{0pt}
\setlength\fboxrule{0pt}
% trim=l b r t
\begin{figure*}[t]
\label{fig:utility}
\centering
\subfloat[]{\label{fig:1a}\fbox{\includegraphics[width=\figwid]{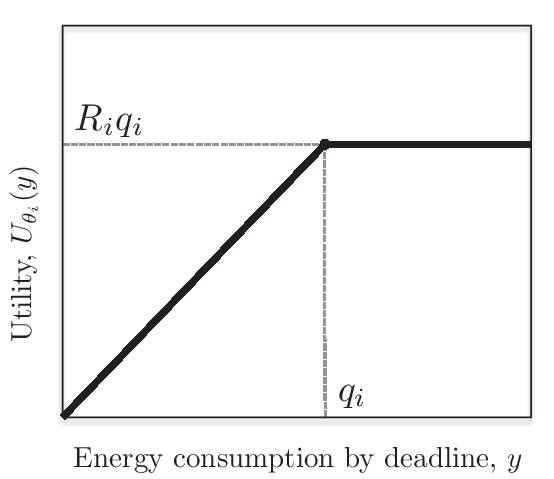}}}\quad
\subfloat[]{\label{fig:1b}\fbox{\includegraphics[width=\figwid]{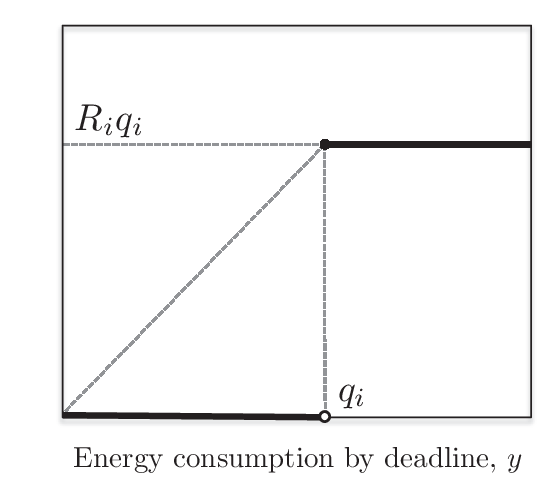}}}\quad
\subfloat[]{\label{fig:1c}\fbox{\includegraphics[width=\figwid]{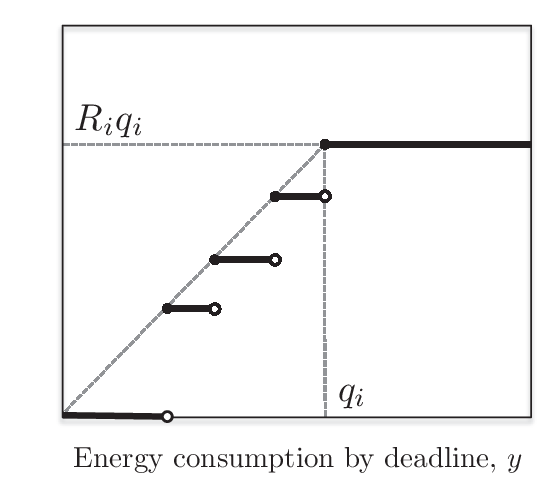}}}\quad
\subfloat[]{\label{fig:1d}\fbox{\includegraphics[width=\figwid]{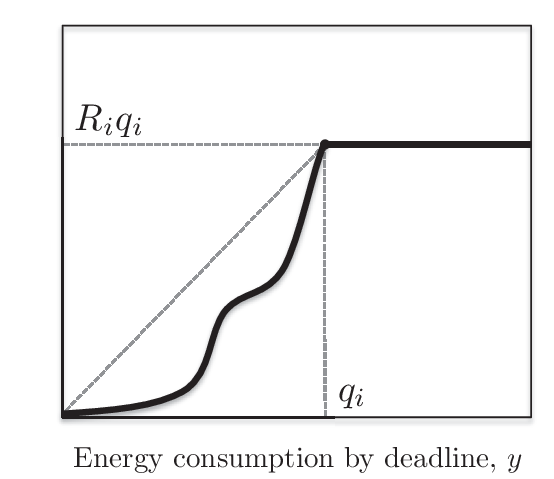}}}
\caption{Four examples of utility functions satisfying Assumption \ref{a:convex}.}
\label{fig:1}
\vspace{-.2in}
\end{figure*}

We consider a model involving a continuum of infinitesimal consumers, indexed by $i \in [0,1]$. Since each individual consumer's action has no influence
on the price, she will act as a price taker. Such an assumption is
reasonable and commonly employed in the context of retail electricity markets \cite{Chao87, Oren90, SO93, Varaiya92}, where it is not uncommon for an electric power utility to service $10^4$ to $10^6$ customers -- a setting in which each consumer herself is too small to influence the price.

%
% This
%setting is reasonable  in the context of retail electricity markets, as
% each consumer herself is too small to
%influence the price. We note that the price-taking assumption is
%somewhat standard in the literature on mechanism design for queueing
%systems, where each consumer is  assumed to be either
%infinitesimally small \cite{MW90,M00,HX09,A13} or price-taking
%\cite{LL97}.

We propose a consumer utility model yielding a  preference ordering on
deadlines, where the longer the delay in consumption, the smaller the utility derived from consumption.
%
%\red{It is natural to assume that the \emph{utility} derived
%from electricity consumption is non-increasing in the delivery
%deadline $k$. Namely, the longer the delay in consumption, the smaller the utility derived from consumption.}
%
 In particular, we assume that each consumer has a
single deadline preference. More precisely, a consumer
with deadline preference $k$ incurs no loss of utility by deferring
 consumption until deadline $k$ and derives zero utility
for any consumption thereafter. This assumption is reasonable for
electric loads such as plug-in electric vehicles (PEVs), dish washers,
and laundry machines, as customers commonly require only that such loads fully execute before a specific time.
Such intuition lends itself to the following definition of \emph{consumer type}.

\begin{definition}[Consumer type]
The \emph{type of consumer} $i$ is defined as a triple
\vspace{-.03in}
\begin{align*}
\theta_i = (k_i,R_i,q_i),
\end{align*}
which consists of her \emph{deadline $k_i \in \Nset$, marginal
utility $R_i \in \Rset_+$}, and \emph{maximum demand} $q_i \in \Rset_+$.
\end{definition}
\noindent Consumer
$i$'s utility function depends only on her type $\theta_i$ and is
assumed to satisfy the following conditions.

\begin{assumption}\label{a:convex}
A consumer  of type $\theta_i$ derives a utility that depends only on
the total energy consumed by her true
deadline $k_i$. The  \emph{utility function} $U_{\theta_i}: \Rset_+ \rightarrow \Rset_+$ is assumed to be non-negative and
non-decreasing over $[0,q_i]$ with
$$
U_{\theta_i}(y) \le yR_i, \quad  \text{for all}  \ \ y \in [0,q_i],
$$
where $R_i= U_{\theta_i}(q_i)/q_i$. We also assume that $U_{\theta_i}(y)=U_{\theta_i}(q_i)$ for all $y  \ge q_i.$  $\hfill \square$
\end{assumption}
Note that the marginal utility $R_i$ associated with a consumer type $\theta_i$ is defined as the ratio
of the maximum utility $U_{\theta_i}(q_i)$ to the maximum demand $q_i$.
It is also worth mentioning that Assumption \ref{a:convex} prevents the treatment of general concave utility functions, as incentive compatibility may fail to hold for certain concave utility functions. We refer the reader to Remark \ref{r:ca} for a more detailed discussion.

\begin{example}[Consumer utility functions]
It is worth emphasizing that Assumption \ref{a:convex} accommodates a large family of utility functions. Indeed, every utility function that lies below the piecewise affine function $U_{\theta_i}(y)=  R_i \min\{y,q_i\}$
satisfies Assumption \ref{a:convex}. Figure 1 depicts several such utility functions.
The utility function depicted in Figure 1(b) characterizes electric loads with  `all or nothing' utility characteristics.  For example, in the context of electric vehicle charging, $q_i>0$ can be interpreted as the minimum amount of battery charge required by the consumer in order to fulfill her next trip. Naturally, the
consumer obtains zero utility if the battery level is below this threshold.
Many job-oriented appliances, such as dishwashers, are similarly modeled.
The utility function in Figure 1(c) is a natural generalization of that in 1(b), and could
 capture the utility function of a consumer having  multiple all-or-nothing  jobs requiring completion before a common deadline.
 Finally, the utility function in Figure 1(d)
is a general nonlinear utility function that satisfies Assumption
 \ref{a:convex}. $\hfill \Box$
\end{example}

We let $\Theta$ denote the \emph{set of all possible consumer
types}, which is assumed to be finite. Let $\rho: \Theta \rightarrow
[0,1]$ \ denote the \emph{distribution of consumer types} over the
space $\Theta$. In other words, for every $\theta \in \Theta$, there is a
$\rho(\theta)$ fraction of consumers of type $\theta$. It follows
that $\sum_{\theta \in \Theta} \rho(\theta) = 1$.

\begin{definition}[Consumer action]
The \emph{action} of a consumer is a vector
$
\v{a}=(a_{1},\ldots,a_{N})^{\top} \in \Rset_+^N,
 $
where $a_{k}$ denotes the amount of energy that is guaranteed delivery by
deadline $k$. The maximum amount of energy any consumer can
request is \ $Q=\max_{(k,R,q)\in \Theta} \{q\}$. Hence, each
consumer's \emph{action space} is restricted to $ \Acal = \left\{
\v{a} \in \Rset_+^{N} \ | \
 \sum_k a_k \leq Q \right\}$.
\end{definition}

It follows from the above definition  that $q \leq Q$ for
 every type $\theta=(k,R,q)$. In other words, it is feasible for every consumer to request
 her maximum demand $q$.
 Given a fixed scheduling policy and pricing scheme, a consumer's {\it strategy} $\varphi: \Theta \to
\Acal$ maps her type into an action.
 In the \yunjian{following} analysis, we will be concerned with
identifying conditions on both the scheduling policy and pricing
scheme that lead to efficient allocations, while inducing
consumers to truthfully reveal their underlying deadline
preferences. A consumer of type $\theta = (k, R, q)$  is defined to be \emph{truth-telling} if she requests $q$
units of energy at her true deadline $k$, and
nothing else. We make this notion precise in the following definition.

\begin{definition} [Truth-telling]\label{def:truth}
A consumer of type $\theta = (k, R, q)$  is defined to be \emph{truth-telling} if her strategy $\v{a}^*=\varphi^*(\theta)$  satisfies $a_j^* = 0$ for  all $j \neq k$ and  $a_k^* = q$.
%$a_j^* = q \cdot \v{1}_{\{ j =k\}}$
%\begin{equation}\label{eq:truth}
%\displaystyle a_j^* = \left\{ \begin{array}{ll} q, & \ \
%j = k, \\[5pt]
%\displaystyle 0, & \ \ j \neq k.
%\end{array} \right.
%\end{equation}
%\begin{align} \label{eq:truth}
%a_j^* = q \cdot \v{1}_{\{ j =k\}}
%\end{align}
%for all $j=1, \dots, N$.
\end{definition}

\yunjian{Indeed, Definition \ref{def:truth} is equivalent to the standard definition of truth-telling in  which a consumer is required to report her exact type.
In Corollary \ref{coro:eff}, we show that social welfare is maximized at a competitive equilibrium, even if the supplier does not have knowledge of  the consumers' utility functions (including the parameter $R$), as long as $R$ is no less than
  the marginal cost of firm supply $c_0$ (cf. Assumption \ref{A:IC}). In other words, a consumer  of type $\theta = (k,R,Q)$ only needs to report her true deadline $k$ and maximum demand $q$ to the supplier, and the profit maximization problem of a (price-taking) supplier naturally yields social optimality.
  In Definition \ref{def:truth}, we require that the consumer requests the entirety of her demand $q$ to be delivered by her true deadline $k$. This is without loss of generality,  because the consumer's utility depends only on the total energy consumed by her true
deadline $k$ (cf. Assumption \ref{a:convex}).}

It is worth mentioning that under an arbitrary scheduling policy and pricing scheme,
it is indeed possible that a consumer's best response is untruthful \yunjian{(cf. Remark \ref{rem:assump} and the discussion following Definition \ref{Def:ICM})}. Our aim is to provide an explicit characterization of a scheduling policy and pricing scheme that implement truth-telling as a dominant strategy for every consumer $i$.

Given the
collection of consumer types $\v{\theta}=\{\theta_i\}_{i \in [0,1]}$ and
a strategy profile $\v{\varphi}=\{\varphi_i\}_{i \in [0,1]}$, the
aggregate demand bundle $\bx$ is given by the mapping
\begin{equation}\label{eq:demand}
\bx=\v{d}(\v{\theta},\v{\varphi})= \int_{i \in [0,1]}
\varphi_i(\theta_i) \eta(di) ,
\end{equation}
where $\eta$ is the Lebesgue measure defined over $[0,1]$, and
$\v{d}=(d_1,\ldots,d_N)$ maps $(\v{\theta},\v{\varphi})$
into an {$N$}-dimensional nonnegative vector.\footnote{Note that we have
implicitly assumed that for every $\v{\theta}$, the
function $\v{\varphi} = \{\varphi_i(\theta_i)\}_{i \in [0,1]}$ is
Lebesgue integrable in $i$. This assumption holds, for example, under a
symmetric strategy profile according to which all consumers of the
same type take the same action.}
Under the \emph{truth-telling strategy profile}
$\v{\varphi}^*=\{\varphi^*_i\}_{i \in [0,1]}$  specified in Definition
\ref{def:truth}, the aggregate demand bundle
 simplifies to
\begin{equation}\label{eq:demandB}
x_j^* =d_{j}(\v{\theta},\v{\varphi}^*)= \sum_{\theta \in \Theta}
q \cdot \rho(\theta)   \cdot \mathbf{1}_{\{j =k\}}
\end{equation}
for all $j=1, \dots, N$, where $\theta=(k,R,q)$. Here, $\mathbf{1}_{\{\cdot\}}$ denotes the indicator function for the event in the subscript.

\old{
\begin{remark}
For the remainder of the paper, we restrict our attention to
aggregate demand bundles $\mathbf{x}$ resulting from truth-telling
consumer populations, where $\mathcal{X}$ denotes the set of all
such demand bundles -- a finite subset of $\Rp^{K+1}$.
\end{remark}}

\subsection{Consumer Surplus}
We are now in a position to characterize the expected surplus
derived by a consumer. It depends on:  (i) her own type and strategy,
(ii) the remaining consumers' types and strategy profile, (iii) the pricing scheme, and of
course, (iv) the scheduling policy employed by the supplier.
%Before formally characterizing an individual consumer's surplus in our
%model,
Before proceeding, we require the definition of pertinent notation.
We define the random
variable \ $\omega^\pi_{k,i}(\bx,\v{a})$  \ to denote the amount of energy delivered to consumer $i$ by
stage $k$ given her requested bundle $\v{a}$ and  aggregate
demand bundle $\bx$. 	This random variable,
which naturally depends on the scheduling policy $\pi$ employed by the
supplier, is formally defined in  Appendix \ref{app:intra}.\footnote{Note that we have allowed the supply
$\omega^\pi_{k,i}(\bx,\v{a})$ to depend explicitly on the consumer index
$i$, as the supplier may employ a scheduling policy that depends on the consumer
index.}
%This dependency on the scheduling policy is made
%precise in Section \ref{sec:supply_model}.

We require that the
requested quantities are always supplied by their corresponding
deadlines and the total quantity delivered to a consumer never exceeds said consumer's total demand.
More formally, we require for each consumer $i\in [0,1]$ taking action $\v{a} \in \Acal$ that
\begin{equation}\label{eq:commit}
\sum\nolimits_{t=1}^k a_{t} \ \leq  \ \omega^\pi_{k,i}(\bx,\v{a}) \ \le \ \sum\nolimits_{t=1}^N a_{t}
\end{equation}
%$$    \omega_{k,i}(\bx,\v{a}') \in \left[ \sum\nolimits_{t=0}^k a'_{t} \ , \ \ \sum\nolimits_{t=0}^K a'_{t} \right] \quad \text{almost surely,}$$
with probability one, for all  aggregate demand bundles $\bx$ and $k = 1,\dots,N$.

%\red{With market efficiency considerations in mind, it is important to understand when a consumer has incentive to misreport her underlying deadline preference.
%Consider a consumer $i$ of type $\theta_i  = (k_i,R_i,q_i)$ facing
%a particular scheduling policy and pricing scheme that satisfy all prior assumptions. Because of monotonicity of prices
%and the service guarantee provided by \eqref{eq:commit}, said
%consumer $i$ has no incentive to request any quantity of energy before her true
%deadline $k_i$. However, if the price of energy associated with later deadlines in low enough, said consumer may have an incentive to report
%a false later deadline if early delivery  is likely (i.e. with high probability) under the specified scheduling policy.
%Intuitively, a consumer $i$ might have incentive to report a false bundle if the reduction in
%total expenditure derived by requesting energy with later
%deadlines  exceeds the expected
%loss of utility incurred by a shortfall in the amount of energy delivered
%by her true deadline $k_i$. Thus, it is essential that, in jointly designing a pricing scheme and scheduling policy, the price differential across deadlines correctly balance the likelihood of early supply under the specified scheduling policy.}

In order to formally define
and analyze incentive compatibility of our proposed market mechanism, we
 require a definition of the expected surplus derived by a consumer under a particular strategy. Given a scheduling policy $\pi$ and a pricing scheme $\v{\kappa}$ \yunjian{(that maps every aggregate demand bundle $\bx$
into a menu of deadline differentiated prices \v{p})}
 employed by the supplier,
and all consumer types $\v{\theta}=\{\theta_i\}_{i \in [0,1]}$,
consumer $i$ receives an \emph{expected surplus} (payoff) under a strategy
profile $\v{\varphi}=\{\varphi_i\}_{i \in [0,1]}$ that is given by
\begin{equation}\label{eq:VkM}
v^\pi_{i}(\v{\theta},\v{\varphi},\v{\kappa}) \  = \ \E \left\{
U_{\theta_i}\Big(\omega^\pi_{k_i,i}(\bx,\v{a})\Big)\right\} \ - \
 \v{\kappa}(\bx)^{\top} \v{a}.
\end{equation}
Here, $k_i$ is consumer $i$'s true deadline, $\bx=\v{d}(\v{\theta},\v{\varphi})$  is the aggregate demand
bundle,
$\v{a}=\varphi_i(\theta_i)$ is the action taken by consumer $i$, and $\v{\kappa}(\bx)$ is
the price bundle set according to the pricing scheme $\v{\kappa}$ at the aggregate demand $\bx$.
Expectation is taken with respect to
the random variable $\omega^\pi_{k_i,i}(\bx,\v{a})$.

%Since each individual consumer is assumed to be infinitesimal in size, and thus have no influence on the price,

Clearly, a scheduling policy $\pi$ together with a pricing scheme $\v{\kappa}$
 defines a game for the
 %(price-taking)
consumer population,
with each individual consumer's payoff expressed in \eqref{eq:VkM}.
We note that this is an {\it aggregative game}, where the payoff
function of each player depends on the population's strategy profile only through the sum of their actions -- the aggregate demand bundle $\v{x} = \v{d}(\v{\theta},\v{\varphi})$.
We can therefore rewrite consumer $i$'s payoff in \eqref{eq:VkM} in a form that depends on other players' strategies only through their sum $\v{x}$. Namely,
\begin{equation}\label{eq:pi}
V^\pi_{i}(\theta_i,\varphi_i,\v{x},\v{\kappa})=\ \E \left\{
U_{\theta_i}\Big(\omega^\pi_{k_i,i}(\bx,\v{a})\Big)\right\} \ - \
 \v{\kappa}(\bx)^{\top} \v{a},
\end{equation}
where
$\v{a}=\varphi_i(\theta_i)$.
 Moreover, it follows from the \blue{service constraint in \eqref{eq:commit}} that under
the truth-telling strategy $\varphi^*_i$, the payoff derived by
consumer $i$ simplifies to the deterministic quantity
\begin{equation}\label{eq:VkMt}
V^\pi_{i}({\theta}_i,\varphi^*_i, \bx,\v{\kappa}) \ = \
U_{\theta_i}(q_i) \ - \  \kappa_{k_i}(\bx) q_i,
\end{equation}
where $\theta_i=(k_i,R_i,q_i)$.
%, and $\bx$ is the aggregate demand of all consumers other than $i$.
 It is important to note that the
expression in \eqref{eq:VkMt} does not depend on the types and
strategies of the other consumers, or the scheduling policy used by
the supplier, as long as the service constraint in \eqref{eq:commit} is respected.

Bayesian Nash equilibrium may not be a plausible solution
concept to explore for this game, as it requires each individual consumer to
have information regarding the distribution of other consumers' types, as well as knowledge of the probability
distribution of $\omega^\pi_{k_i,i}(\bx,\v{a})$, which in turn depends on the
distribution of the intermittent supply process $\v{s}$. We circumvent such informational assumptions by
focusing our analysis around a  stronger solution concept -- namely, the {dominant strategy
equilibrium} of the game. A strategy $\varphi_i$ is a \emph{dominant
strategy} for consumer $i$ of type $\theta_i$, if it maximizes her
expected payoff regardless of the actions taken by the other
consumers. We have the following definition.

\begin{definition}[Dominant strategy] \label{Def:DS}
A strategy
$\varphi_i$ is a \emph{dominant
strategy} for a consumer $i$ of type $\theta_i$ if
$$
V^\pi_{i}(\theta_i,\varphi_i,\v{x},\v{\kappa}) \ \ge \
V^\pi_{i}(\theta_i,\varphi'_i,\v{x},\v{\kappa}), \quad  \forall  \
\varphi'_i, \quad \forall \ \bx \in \Rset_+^{N}.
$$
\end{definition}

%\begin{definition}[Dominant strategy]
% A strategy $\varphi_i$ is said to be a \textbf{dominant strategy} for consumer $i$ of type $\theta_i$, if it maximizes her
%expected payoff regardless of the actions taken by the other
%consumers, i.e.
%$$
%V^\pi_{i}(\theta_i,\varphi_i,\v{x},\v{\kappa}) \ \ge \
%V^\pi_{i}(\theta_i,\varphi'_i,\v{x},\v{\kappa}), \qquad  \forall  \
%\varphi'_i, \qquad \forall \ \bx \in \Rset_+^{N}.
%$$
%\end{definition}

 We note that in Definition \ref{Def:DS}, a change of an individual consumer $i$'s strategy
has no influence on the aggregate demand bundle $\bx$, because each consumer is assumed to be infinitesimal in size, and the aggregate demand bundle is given by Eq. \eqref{eq:demand}.
In Definition \ref{Def:ICM}, we define a  mechanism to be \emph{incentive compatible} if it implements truth-telling in dominant strategies.
Surprisingly,  we show in Section \ref{sec:res2} that a mechanism consisting  of an earliest-deadline-first (EDF) scheduling policy in combination with a marginal cost pricing scheme
is incentive compatible and induces an efficient competitive equilibrium between the supplier and consumers.

\old{
\begin{remark}[An alternative formulation]
One can define an alternative game among the consumer population by leveraging on an additional assumption that the supplier has a priori knowledge of the aggregate demand curve, but not necessarily the type of each individual consumer.
We note that this is the typical setting considered in the literature on priority
pricing for queuing systems \cite{MW90,RB95,LL97,M00}.
In this setting, the supplier can simply announce a predetermined price bundle $\v{p} \in \Pcal$ (as opposed to a pricing scheme $\v{\kappa}(\cdot)$), which also defines a game between the consumers. Our results on incentive compatibility outlined in Section \ref{sec:res2} imply that with the knowledge of the aggregate demand curve, the supplier can simply set the price bundle to coincide with the intersection of the marginal cost supply and demand curves, thus inducing
consumers' truth-telling behavior at an efficient competitive equilibrium between the supplier and consumers. We discuss this alternative approach in more detail in Section \ref{sec:CE}.
\end{remark}}

%%%%%%%%%%%%%%%%%%%%%%%%%%%%%%%%%%%%%%%%%%%%%%%%%%%%%%%%%%%%%%%%%%%%%%%%%

%\vspace{-3mm}

\section{Model of Supply} \label{sec:supply_model}

 As one of the primary thrusts of this paper is the characterization of a competitive equilibrium between the supplier and consumers and its efficiency properties, we now consider the behavior of a \emph{price-taking supplier}, whose aim is to maximize his expected profit given a predetermined price bundle. The expected profit derived by a supplier equals the revenue derived from the sale of a bundle of deadline differentiated energy quantities less the expected cost of firm supply required to service said bundle. In determining  his supply curve under price taking behavior, the supplier's objectives are two-fold:

\begin{itemize} \setlength{\itemsep}{.06in}
\item  \emph{Scheduling.} Determine an \emph{optimal scheduling policy}  to causally allocate the intermittent supply across the deadline differentiated consumer classes, in order to minimize
the expected cost of firm supply required to ensure that demand
bundles are served by their respective deadlines.
%\blue{Yunjian: the following sentence needs work?} Naturally, the scheduling policy will be parameterized by the aggregate demand bundle.

\item  \emph{Pricing.} Given the optimal scheduling policy, determine a \emph{marginal cost supply curve} that specifies the bundle of energy he is willing to supply at every price bundle.
\end{itemize}

Essentially, the specification  of the supplier's marginal cost curve requires the explicit characterization  of the gradient of the supplier's minimum expected cost of firm supply. This, in turn, requires the specification of supplier's optimal scheduling policy, which is proven to have an earliest-deadline-first (EDF) structure in Theorem \ref{thm:opt_pol_app}. 
With the optimal scheduling policy in hand, we derive a closed-form expression for the supplier's  marginal cost curve  in Theorem \ref{thm:opt_sup}.

%With these, fuwhich is used to construct  an efficient
%competitive equilibrium between supply and demand.
%
%
%requires the explicit characterication of the suppliers 
% amounts to the explicit solution of a two-stage stochastic program,
% whose expected recourse cost is the optimal value of a constrained scheduling problem
% parameterized by the aggregate demand bundle.
%We provide an explicit characterization of such solution 
%
%. The resulting marginal cost supply curve 
%
%The 

\subsection{Feasible Scheduling Policies}\label{sec:fs}
We now offer a precise formulation of the supplier's class of feasible scheduling policies given an aggregate demand.
%The interested reader is referred to Appendix \ref{app:sched} for an expanded formulation and technical proofs. 
When considering the problem of scheduling, it is important  to distinguish
between intra-class and inter-class
scheduling. An \emph{inter-class scheduling policy} (denoted by $\sigma$)
represents a sequence of scheduling decisions, which causally allocate the intermittent and firm supply across the deadline-differentiated demand classes. 
%We denote by $\Sigma(\mathbf{x})$ the space of all \emph{feasible inter-class scheduling policies} available for use by the supplier, given an aggregate demand bundle $\bx$.
In addition to inter-class scheduling, the supplier must also determine as to how the available supply is allocated between customers within a given demand  class. As such,
we let the \emph{intra-class scheduling policy} (denoted by $\phi$) represent
 a sequence of scheduling decisions, which dictates how the intermittent and firm supply  made available to a  given demand class is allocated across customers within said class.
% 
% We denote by $\Phi(\sigma)$ the space of all \emph{feasible intra-class scheduling policies} available for use by the supplier,
%which naturally depends on the inter-class scheduling policy $\sigma$ being used.
Finally, we let  $\pi = (\sigma, \phi)$
%$$\pi = (\sigma, \phi), \quad \text{where} \ \   (\sigma, \phi) \in \Sigma(\mathbf{x}) \times \Phi(\sigma)$$ 
denote  the \emph{joint inter-class and intra-class scheduling policy} employed by the supplier. In what follows, we formally define the space of feasible inter-class and intra-class scheduling policies.
%A formal definition of feasible intra-class policies is developed in Appendix B of \cite{BX15}.

\

\subsubsection{Inter-Class Scheduling Policies}  \label{sec:inter}
We  characterize the optimal inter-class scheduling policy as a
solution to a constrained stochastic optimal control problem. First,
we define the system \emph{state} at period $k$ as the pair
$(\mathbf{z}_k, s_k) \in \Rset^{N}_+ \times \Rset_+$, where the
vector $\mathbf{z}_k$ denotes the residual demand requirement of the
original aggregate demand bundle $\bx$ after having been serviced in
previous periods $0,1,\cdots,k-1$.
Define as the \emph{control input} the vectors $\mathbf{u}_k, \mathbf{v}_k \in \Rset_+^{N}$,
which denote (element-wise in $j$) the amount of intermittent and firm supply allocated to demand class $j$ at period $k$, respectively. Naturally then, the state of residual demand evolves according to the discrete time \emph{state equation}:
\begin{align}
\mathbf{z}_{k+1} \ = \ \mathbf{z}_{k} - \mathbf{u}_{k} - \mathbf{v}_{k}, \qquad k=0,\dots, N-1,
\end{align}
 where the process is initialized with \ $\v{z}_0 = \v{x}$.
The delivery deadline constraints manifest in a sequence of nested constraint sets $ \Rset_+^N \supset \Zcal_1 \supseteq \Zcal_2 \supseteq \cdots \supseteq \Zcal_N = 0$ converging to the the origin, where the set $\Zcal_k$ characterizes the \emph{feasible state space} at stage $k$.
 More precisely,
\[
\Zcal_k \  =  \ \{ \v{z} \in \Rset_+^N  \ | \  z^j =0, \ \ \forall \
j \leq k
\}.
 \]
 In other words, the feasible state space is such that each
demand class is  fully serviced by its corresponding deadline. We
define as the \emph{feasible input space} at stage $k$ the set of
all inputs belonging to the set
\[ 
\Ucal_k(\v{z}, s)  =   \{ (\v{u}, \v{v})  \ | \ \mathbf{1}^{\top} \v{u} \leq s \  \text{and}  \ \v{z} - \v{u} - \v{v} \in \Zcal_{k+1} \}, 
\]
which ensures one-step state feasibility and that the total allocation of renewable supply does not exceed availability at the current stage. In characterizing the feasible set of causal scheduling policies, we restrict our attention to those policies with \emph{Markovian information structure}, as opposed to allowing the control to depend on the entire history.
This is without loss of optimality, \old{as we will later show in Lemma \ref{lemma:oracle}, which reveals, more strongly,}  since Markovian policies are capable of performing as well as the optimal oracle policy.  We describe the scheduling decision at each stage $k$ by the functions
\[ \v{u}_k = \v{\mu}_k(\v{z}, s) \quad  \text{and} \quad \v{v}_k = \v{\nu}_k(\v{z}, s), \]
where $\v{\mu}_k: \Zcal_k \times \Scal \rightarrow \Rset_+^N$   and $\v{\nu}_k :\Zcal_k \times \Scal \rightarrow \Rset_+^N$.
A \emph{feasible inter-class scheduling policy} is any finite sequence of scheduling decision functions  $\sigma = (\v{\mu}_0, \dots, \v{\mu}_{N-1}, \v{\nu}_0, \dots, \v{\nu}_{N-1})$ such that 
\begin{align*}
(\v{\mu}_k, \v{\nu}_k)(\v{z}, s) \in \Ucal_k(\v{z}, s), \quad \forall  \ \ (\v{z}, s) \in \Zcal_k \times \Scal
\end{align*}
and  time periods $k=0,\dots, N-1$.  We will occasionally write the state and control process as \ $\{\v{z}_k^{\sigma}\}, \ \{\v{u}_k^{\sigma}\},$ \ and  $ \{ \v{v}_k^{\sigma}\}$  to emphasize their dependence on the inter-class policy $\sigma$, unless otherwise clear from the context.

Given an aggregate demand bundle $\bx$, we denote by $\Sigma(\mathbf{x})$ the \emph{space of all feasible inter-class scheduling policies} available for use by the supplier.

%
%An important assumption we will make on the supply side is the absence of an upper bound on the amount of energy
%the supplier can deliver to a customer within any given time period.
%Although appearing strong at first glance, such an assumption is not far from
%practice, as batteries with high power to energy ratios are rapidly becoming available for electric vehicles. For example,
%the lithium-ion titanate batteries are capable
%of recharging  to 95\% of full capacity within approximately ten minutes \cite{BS09,EA10}. We also note that
%fast (DC) chargers that can fully charge an electric vehicle within half an hour are being installed
%in a variety of public locations including
%parking lots, shopping centers, hotels, theaters, and restaurants \cite{YK13}.

\

\subsubsection{Intra-Class Scheduling Policies}
\label{sec:intra} 
Recall that an intra-class scheduling policy $\phi$
determines the allocation of  supply to specific consumers within each
deadline-differentiated demand class, where the supply available to
each demand class is determined by the inter-class policy.  Given an inter-class scheduling policy $\sigma \in \Sigma(\mathbf{x})$,  we denote by $\Phi(\sigma)$ the \emph{space of all feasible
intra-class scheduling policies}.\footnote{As the formal definition of feasible intra-class policies is relevant only to the derivation of technical proofs, we defer their precise specification to Appendix \ref{app:intra}  in order to maintain continuity in exposition.}

%We denote by $\pi = (\sigma, \phi)$ the joint inter and intra-class scheduling policy employed by the supplier.

It is important to note that the  supplier's expected
profit depends only on the inter-class scheduling policy being used, and is 
invariant under the family of feasible intra-class scheduling policies.
This follows from the supplier's cost indifference to supply allocation
between consumers within a given demand class.
Therefore, in
characterizing the optimal scheduling policy  for the supplier, we
restrict our attention to the characterization of optimal inter-class policies. 
The distinction between inter-class and intra-class scheduling is, however, important, as the choice of intra-class scheduling policy $\phi$ can affect the
probability distribution of each consumer's random supply $\omega_{k_i,i}^{\pi}(\bx,\v{a})$. This can, in turn, influence consumer purchase decisions.

\subsection{Optimal Scheduling Policy}\label{sec:os}

We define the \emph{expected profit} $J(\v{x},  \v{p}, \sigma)$ derived by a supplier as the revenue derived from an aggregate demand bundle $\v{x}$ less the expected cost of servicing said demand bundle under a feasible inter-class scheduling policy $\sigma \in \Sigma(\v{x})$. More precisely, let
\begin{eqnarray*} 
J(\v{x},  \v{p}, \sigma) \ = \ \v{p}^{\top} \v{x} \ - \ Q(\v{x}, \sigma),
\end{eqnarray*}
where $Q$ denotes the expected cost of firm generation incurred servicing $\v{x}$ under a feasible policy $\sigma \in \Sigma$.
It follows that
\begin{align}
Q(\v{x}, \sigma)  \ = \  \sum_{k=0}^{N-1} \mathbb{E}  \left\{ \v{c}^{\top} \v{v}_k^{\sigma}  \right\},
\end{align}
where $\v{c} = (c_0,\dots, c_0)$. We wish to characterize inter-class scheduling policies, which lead to a minimal expected cost of firm supply. Accordingly, we have the following definition of optimality.
\begin{definition}[Optimal policy]
The inter-class scheduling policy $\sigma^* \in \Sigma(\v{x})$ is defined to be \emph{optimal} if
\begin{align} \label{opt_sched11}
Q(\v{x}, \sigma^*) \ \leq \ Q(\v{x}, \sigma),  \quad \forall \ \ \sigma \in \Sigma(\v{x}).
\end{align}
We denote by $Q^*(\v{x}) = Q(\v{x}, \sigma^*)$ the \emph{minimum expected cost of firm supply.}
\end{definition}

The following result provides an explicit characterization of an optimal inter-class scheduling policy.

\begin{theorem}[Earliest-Deadline-First] \label{thm:opt_pol_app}
Given an aggregate demand bundle $\v{x} \in \Rset_+^N$, the optimal scheduling policy $\sigma^* \in \Sigma(\v{x})$ is given by:
$$
\begin{array}{l}
\mu_k^{j,*}(\v{z}, s)  =  \min\left\{ z^j ,   s - \sum\nolimits_{i=1}^{j-1} \mu_k^{i,*}(\v{z}, s)  \right\}, \\[8pt]
\nu_k^{j,*}(\v{z}, s)=(z^j - \mu_k^{j,*}(\v{z}, s))\cdot \mathbf{1}_{\{k=j-1\}},
\end{array}
$$
for  $j = 1,\dots, N$, $k = 0,\dots, N-1$,   and $(\v{z}, s) \in \Zcal_k \times \Scal$.
\end{theorem}

Theorem \ref{thm:opt_pol_app} is intuitive.\footnote{We omit a formal proof of  Theorem  \ref{thm:opt_pol_app}, given its immediacy in derivation upon examination of scheduling problem's dynamic programming equations.}
The optimal inter-class policy $\sigma^* \in \Sigma(\v{x})$ is such that  the intermittent supply $s_k$ available at each period $k$ is allocated to those unsatisfied demand classes with \emph{earliest-deadline-first} (EDF), while the firm supply is allocated to a demand class only when the EDF allocation of intermittent supply is insufficient to ensure its deadline satisfaction. In other words,  the firm supply is used only as a last resort. 
For the remainder of the paper, we refer to $\sigma^*$ as the \emph{EDF scheduling policy.} 

\begin{remark}[Implementation requirements]
An attractive property of the EDF  scheduling policy  $\sigma^*$ is that it is \emph{distribution-free}. Namely, it can be implemented by the supplier without requiring explicit knowledge of the underlying probability distribution of the intermittent supply process. On the other hand, a  potential limitation  of the proposed scheduling policy is the centralization of information exchange it entails, as it's implementation requires periodic communication between the supplier and each consumer-specific device. Specifically, at every time  period $k$, the supplier must transmit a control signal to each device specifying its energy consumption level in that period. These control signals are determined according to the intra-class scheduling policy $\phi \in \Phi(\sigma^*)$ being used by the supplier.
In return, each device must transmit a measurement signal to the supplier summarizing its state (e.g., residual energy requirement). Naturally, the time scale at which such control schemes can be implemented will 
depend on the number of devices being coordinated and the  available bandwidth of the underlying communication network being used. 
\end{remark}

%We omit a formal proof of  Theorem  \ref{thm:opt_pol_app}, given its immediacy in derivation. Essentially, the crux of the proof centers on showing that the EDF scheduling policy performs almost surely as well as any \emph{non-causal} policy with perfect foresight -- the so called oracle optimal policy.  Technically, the proof relies on showing that any oracle optimal schedule  can be inductively mapped to the EDF schedule without  incurring an increase in the corresponding cost of firm supply, almost surely.

\subsection{Marginal Cost Pricing} \label{sec:opt_price}
Given the EDF characterization of the optimal inter-class scheduling
policy in Theorem \ref{thm:opt_pol_app}, we are now in a position to
characterize the supplier's optimal supply curve under price taking
behavior. We define the  \emph{residual
process} induced by the EDF scheduling policy $\sigma^*$ as
\begin{align}\label{eq:xi}
\xi_{k+1}(\v{x},\v{s}) \ = \  \max\{0, \; \xi_{k}(\v{x}, \v{s})\}   \ + \ s_{k} - x_{k+1},
\end{align}
%\begin{align}\label{eq:xi}
%\xi_{k+1}(\v{x},\v{s}) \ = \   \max\{ \xi_{k}(\v{x}, \v{s}), \; 0\} \ + \ s_{k} - x_{k+1},
%\end{align}
for $k =0, \dots, N-1$, where $\xi_0 = 0$. We denote the entire residual process by $\v{\xi} = (\xi_0, \dots, \xi_N)$, omitting its dependency on $\bx$ and $\mathbf{s}$  when it is clear from the context.
A positive residual ($\xi_k
>0$) represents the amount of intermittent supply leftover after
having serviced demand class $k$ by its
deadline $k$, according to the EDF inter-class scheduling policy $\sigma^*$. A negative residual  ($\xi_k \leq 0$) represents the
amount by which the intermittent supply fell short --  or, equivalently,
the amount of firm supply required to ensure satisfaction of the
demand class $k$.
%While $\xi_k(\v{x},\v{s})$
%depends on both the aggregate demand bundle $\bx$ and intermittent supply process
%$\mathbf{s}$, we omit this dependency when it is clear from the context and compactly write \ $\v{\xi} = (\xi_0, \dots, \xi_N) \in \Rset^{N+1}$.
Using this newly defined process, we have the
following characterization of the minimum expected cost of
firm supply under EDF scheduling. First, we require a technical assumption.

\begin{assumption} \label{assum:dist}
The joint probability distribution of the intermittent supply process $\v{s}$ is assumed to be absolutely continuous and have compact support. $\hfill \Box$
\end{assumption}

\begin{lemma} \label{lem:profit}
\blue{Suppose that Assumption \ref{assum:dist} holds.} The minimum expected recourse cost  $Q^*(\v{x})$  derived under an aggregate demand bundle $\v{x} \in \Rset_+^N$ and  EDF scheduling policy  \ $\sigma^* \in \Sigma(\v{x})$   \ satisfies
\begin{align} \label{eq:opt_val}
Q^*(\v{x}) \ = \   \mathbb{E} \left\{c_0 \sum_{k=1}^N  - \min\{0, \; \xi_{k}\}  \right\},
\end{align}
and is convex and differentiable in $\v{x}$ over $[0,\infty)^N$.
\end{lemma}
Lemma \ref{lem:profit} admits a  natural interpretation. Namely, the minimum expected  cost of firm supply is equivalent to the amount by which the intermittent supply is expected to fall short for each demand class under EDF scheduling. Moreover, it follows readily that the expected profit  $J(\v{x}, \v{p}, \sigma^*)$  derived under EDF scheduling is also differentiable and concave in $\v{x}$. As such, any allocation $\v{x}$ satisfying the first order condition,
\begin{align} \label{eq:foc}
\nabla_{\v{x}} \;  J(\v{x}, \v{p}, \sigma^*)^{\top}( \v{x} - \v{y}) \  \geq \ 0  \qquad \text{for all}  \ \ \ \v{y} \in \Rset_+^N,
\end{align}
is profit maximizing for the supplier given a  price bundle $\v{p} \in \Pcal$.  Accordingly, we provide an explicit characterization of the supplier's \emph{marginal cost supply curve} in the following theorem.

%Identifying a price bundle that induces a competitive equilibrium between supply and demand requires, first, a characterization of the supplier's optimal supply curve. Namely, given a price bundle $\v{p}$, the supplier computes
%\begin{align}
%\v{x}^* \ \in \  \arg \max_{\v{x} \in \Rset_+^N}  \  J(\v{x}, \v{p}, \pi^*),
%\end{align}
%where $\v{x}^*$ denotes the supplier's profit maximizing supply allocation -- an explicit function of the price bundle $\v{p}$. By concavity of $J(\v{x}, \v{p}, \pi^*)$, one can readily compute necessary and sufficient conditions for optimality, which are made precise in the following Theorem \ref{thm:opt_sup}.

\begin{theorem}[Marginal cost supply curve] \label{thm:opt_sup} \blue{Suppose that Assumption \ref{assum:dist} holds.} An allocation $\bx$ is profit maximizing for a given price bundle $\mathbf{p}$
if   $$\v{p}  =  \v{\zeta}(\bx),$$  
where \yunjian{the pricing scheme}  $\v{\zeta} : \Rset_+^{N} \rightarrow \mathcal{P}$  satisfies, 
 \begin{align}\label{eq:eP}
\frac{\zeta_k(\bx)}{c_0}    =
 \mathbb{P}(\xi_k \le 0)  + \sum_{t=k+1}^N \mathbb{P}(\xi_k > 0,\ldots,\xi_{t-1}>0, \xi_{t}\le 0)
\end{align}
for $k=1,\ldots,N$.
\end{theorem}

The \emph{marginal cost pricing scheme} $\v{\zeta}$ specified in Equation \eqref{eq:eP} maps every aggregate demand bundle in $\Rset_+^{N}$
into a menu of deadline differentiated prices in $\mathcal{P}$, \yunjian{which reflects the marginal cost of the supplier}.
Moreover, one can readily interpret such pricing scheme as setting the price $p_k$ for energy with deadline $k$  equal to (up to a  constant factor $c_0$)  the probability that firm supply will be required to service the bundle $\v{x}$ at any subsequent time period $t  \geq k-1$ under the optimal inter-class scheduling policy $\sigma^*$. Naturally, the larger the probability of
shortfall, the larger the price.
%
%
%\begin{example} Consider a two period market ($N=2$) having an intermittent supply process $\mathbf{s} = (s_0,s_1)$ that is independent across time. Given a demand bundle $\bx = (x_1, x_2)$, it is straightforward to show that the deadline differentiated prices $\mathbf{p} =(p_1, p_2)$ defined according to marginal cost supply curve \eqref{eq:eP} are related according to
%\begin{align*}
%p_2 = p_1 - c_0\cdot \mathbb{P}(s_0 \leq x_1)\cdot \mathbb{P}(s_1 \geq x_2).
%\end{align*}
% $\hfill \Box$
%\end{example}
%
%
It is
readily verified that the pricing scheme  $\v{p} \ = \ \v{\zeta}(\bx)$
yields, in general, prices that are nonincreasing in the deadline. More precisely,
%\[  c_0 \ \geq \ \zeta_1(\bx) \ \geq \ \zeta_2(\bx) \ \geq \cdots \geq \ \zeta_N(\bx) \]
\[  c_0 \ \geq \ p_1\ \geq \ p_2 \ \geq \cdots \geq \ p_N \]
for all  $\bx \in
\Rset_+^N$. This property of price monotonicity is consistent with our initial
criterion for constructing such a market system. Namely, the longer a
customer is willing to defer her consumption in time, the less she
is required to pay per unit of energy. And the price of deferrable energy is no greater than the nominal flat rate for electricity, $c_0$.

\section{Incentive Compatibility and Efficiency} \label{sec:res2}

We now establish several important properties of the  market mechanism developed in Section \ref{sec:supply_model}.
In particular, under an additional mild assumption on each consumer's marginal valuation on energy, we show in Sections \ref{sec:IC} and \ref{sec:imp} that a mechanism consisting of earliest-deadline-first (EDF) scheduling in combination with the marginal cost pricing scheme in \eqref{eq:eP} is indeed \emph{incentive compatible},
and \emph{achieves social optimality} at a competitive market equilibrium, respectively.

\subsection{Incentive Compatibility}\label{sec:IC}

With market efficiency considerations in mind, it is important to understand when a consumer has incentive to misreport her underlying deadline preference.
A  mechanism $(\pi,\v{\kappa})$  consisting of a feasible scheduling policy $\pi = (\sigma, \phi)$ and a pricing scheme {$ \v{\kappa}$}
is said to be \emph{incentive compatible} (IC) for
consumers of a particular type, if it is a dominant-strategy for
consumers of this type to be truth-telling. %We have the following definition.

\begin{definition}[Incentive compatibility] \label{Def:ICM}
 A mechanism $(\pi,\v{\kappa})$ is \emph{incentive compatible} for every consumer $i$ of type $\theta_i$, if
\begin{equation*}\label{eq:IC}
V^\pi_{i}(\theta_i,\varphi^*_i, \bx, \v{\kappa}) \ \ge \
V^\pi_{i}(\theta_i,\varphi_i, \bx, \v{\kappa}), \quad
\forall \; \varphi_i, \;\; \forall \; \bx \in \Rset_+^{N},
\end{equation*}
%for all $\varphi_i$ and $ \bx \in \Rset_+^{N}$.
where  $\varphi^*_i$ is the truth-telling strategy
defined in Definition \ref{def:truth}.
\end{definition}

It is worth mentioning that the above definition of dominant strategy incentive compatibility is strong.
It requires that
a consumer would like to reveal her true type regardless of the types and actions of other consumers and the probability distribution on the intermittent
supply process.
Since the optimal price schedule is non-increasing in the deadline, and
demand is guaranteed to be met before the requested deadline, a
consumer $i$ does not have an incentive to request any quantity of
energy before her true deadline $k_i$.
However, if the price of energy associated with later deadlines \yunjian{is} low enough, said consumer may have an incentive to report
a false later deadline if early delivery  is likely (i.e., with high probability) under the specified scheduling policy.
Intuitively, a consumer $i$ will have incentive to misreport its deadline if the reduction in
total expenditure derived by requesting energy with later
deadlines  exceeds the expected
loss of utility incurred by a shortfall in the amount of energy delivered
by her true deadline $k_i$.
 Surprisingly, we show in Theorem \ref{t:IC} that EDF scheduling in combination with marginal cost pricing precludes this possibility.

%
%However, a consumer $i$ may have incentive to report a false deadline $k' > k_i$, if the prices associated with later deadlines ($>k_i$) are low enough.  In particular, a rational (expected-payoff maximizing) consumer
%will report a false later deadline, if the corresponding saving in energy cost exceeds her expected loss of utility (resulting from the possible
%shortfall in the energy received by her true deadline).
%
%
%\red{With market efficiency considerations in mind, it is important to understand when a consumer has incentive to misreport her underlying deadline preference.
%Consider a consumer $i$ of type $\theta_i  = (k_i,R_i,q_i)$ facing
%a particular scheduling policy and pricing scheme that satisfy all prior assumptions. Because of monotonicity of prices
%and the service guarantee provided by \eqref{eq:commit}, said
%consumer $i$ has no incentive to request any quantity of energy before her true
%deadline $k_i$.
%Thus, it is essential that, in jointly designing a pricing scheme and scheduling policy, the price differential across deadlines correctly balance the likelihood of early supply under the specified scheduling policy.}

For the remainder of the paper, we denote by $(\pi^*, \v{\zeta})$ the market mechanism defined by EDF scheduling and marginal cost pricing. More precisely, the scheduling policy $\pi^* = (\sigma^*, \phi)$ consists of the  EDF inter-class scheduling policy $\sigma^* $ and an arbitrary
 feasible intra-class policy $\phi \in \Phi(\sigma^*)$. And, $\v{\zeta}$ denotes the  marginal cost pricing scheme  in Equation \eqref{eq:eP}.

\begin{theorem}[Incentive compatibility] \label{t:IC}
Suppose that Assumptions \ref{a:convex}-\ref{assum:dist} hold.
The mechanism $(\pi^*, \v{\zeta})$  is
\emph{incentive compatible} for all consumers of a type that satisfies $R
\ge c_0$.
\end{theorem}
%\red{We defer the proof of Theorem \ref{t:IC} to Appendix \ref{sec:A}
%for the ease of exposition.}
We further establish in Corollary \ref{coro:eff}  that the mechanism $(\pi^*, \v{\zeta})$  also maximizes social welfare at a unique market equilibrium between the supplier and consumers, if the condition $R\ge c_0$ holds for
every consumer type $\theta \in \Theta$.

\begin{remark} \label{rem:assump} We note that the requirement $R \ge c_0$ is reasonable for
electricity consumers, as their marginal valuation on electricity
consumption is commonly higher than the nominal flat rate for electricity, which in our model is denoted by $c_0$. On this basis, electricity demand is generally modeled as inelastic
\cite{S02,W02}, especially in the short term \cite{ZP07}.
\end{remark}

\begin{remark} \label{r:ca}
We also note that  incentive compatibility may fail to hold  if certain assumptions regarding a consumer's utility function are violated.
First, one can readily show that incentive compatibility may fail to hold for a consumer $i$ of type $\theta_i = (k_i, R_i, q_i)$,
if the marginal cost of firm supply exceeds her marginal valuation of energy -- namely,  $R_i < c_0$. Second, we note that the result in Theorem \ref{t:IC} fails to hold for arbitrary concave utility functions. This is intuitive. Consider a
consumer $i$ having a highly concave utility function. Because of the large underlying concavity in her utility, said consumer
will prefer to
report a false deadline that is later than her true deadline $k_i$, \
 if she can obtain a fraction of her demand before stage $k_i$ with high probability and at a low price.
%
% (but a
%utility close to the maximum, because of the large underlying concavity of her utility
%function) with high probability, and at a much lower price.
\end{remark}

\subsection{Market Equilibrium and Efficiency}\label{sec:imp}
In this section, we show that mechanism consisting of marginal cost pricing together with EDF scheduling,
results in an efficient market equilibrium at
which social welfare (the sum of aggregate consumer surplus and
supplier profit) is maximized.
First, we make an assumption under which we will operate for the remainder of the paper.

%\red{Recall that  $(\pi^*, \v{\zeta})$ denotes the market mechanism defined by EDF scheduling and marginal cost pricing.
%We now discuss the mechanism implementation and its equilibrium properties.}
%
%%
%%We briefly discuss the implementation of the proposed market for
%%deadline differentiated energy services through both a mechanism design and market equilibrium based approach.
\begin{assumption}\label{A:IC}
We assume that every consumer type $(k,R,q) \in \Theta$ satisfies $R \geq c_0$. $\hfill \Box$
\end{assumption}

%
%\red{
%As a result, the proposed market can be implemented through a classical mechanism-design approach:
%
%\
%
%\noindent \textbf{Step 1}. \ The supplier announces the mechanism \ $(\pi^*, \v{\zeta})$ to the population of customers.
%
%\
%
%\noindent \textbf{Step 2}. \  With  common knowledge of the mechanism,
% every consumer $i$ reports her type $\theta_i$ to the supplier, as it is a dominant strategy to do so.
%
%\
%
%\noindent \textbf{Step 3}. \  The price bundle is set to  $\v{p} = \v{\zeta}(\bx^*)$, where  $\bx^* = \v{d}(\v{\theta},\v{\varphi}^*)$ is the truth-telling
%aggregate demand bundle. And each consumer $i$ receives an expected payoff of
%$$
% \E \left\{
%U_{\theta_i}\Big(\omega^{\pi^*}_{k,i}(\v{x}^*,\v{a}^*)\Big)\right\} \ - \v{p}^{\top} \v{a}^*,
%$$
%where  $\v{a}^* = \varphi_i^*(\theta_i)$ is the truth-telling
%action of consumer $i$ defined in Definition \ref{def:truth}.}

See Remark \ref{rem:assump} for a discussion on Assumption \ref{A:IC}.
It follows from Assumption \ref{A:IC} and  Theorem \ref{t:IC} that
under the mechanism $(\pi^*, \v{\zeta})$, it is a dominant strategy for every consumer to be truth-telling.
The aggregate demand resulting from the true-telling population, which we denote by
$\bx^*$, is given by Equation \eqref{eq:demandB}.
Theorem \ref{thm:opt_sup} shows that it is profit-maximizing for a supplier to
meet the aggregate demand $\v{x}^*$ at the price bundle $\v{p} = \v{\zeta}(\bx^*)$. In what follows,
we will show that the resulting quantity-price pair $(\v{x}^*, \v{\zeta}(\bx^*))$ constitutes a market equilibrium that is social welfare maximizing. We first offer a definition of \emph{market equilibrium}.

\begin{definition}[Market equilibrium]\label{def:CE} Let $(\pi,\v{\kappa})$ be a market mechanism consisting of a scheduling policy $\pi$ and pricing scheme $\v{\kappa}$.
Given the types of all consumers $\v{\theta}$, a quantity-price pair
$(\v{x},\v{p})$ is a \emph{market equilibrium} under the mechanism $(\pi,\v{\kappa})$ if the following two
conditions hold.

\vspace{.03in}
\begin{itemize} \setlength{\itemsep}{.05in}
\item[(i)]  It is a dominant strategy for every consumer to be truth-telling under the mechanism $(\pi,\v{\kappa})$, and the
resulting aggregate demand bundle is $\bx$.

\item[(ii)] The aggregate demand bundle $\bx$ together with the scheduling policy $\pi$ maximize
a price-taking supplier's expected profit at price bundle $\v{p} = \v{\kappa}(\bx)$.
% The pricing bundle $\v{p}$ is determined by the pricing scheme
%$\v{\kappa}$ at the aggregate demand $\v{x}$, i.e., $\v{p} = \v{\kappa}(\bx)$.
%Further, given the price bundle $\v{p}$, the aggregate demand bundle $\v{x}$
%and scheduling policy $\pi = (\sigma, \phi)$ maximize the supplier's
%expected profit (cf. its expression in Eq. \eqref{eq:profit22}).
\end{itemize}
\end{definition}

In Definition \ref{def:CE}, the first condition ensures that under the
mechanism $(\pi,\v{\kappa})$, the aggregate demand  (resulting from a truth-telling  non-atomic population) is $\bx$. The
 second condition
requires that given the price bundle $\v{p}= \v{\kappa}(\bx)$
(which is determined according to the pricing scheme $\v{\kappa}$ and aggregate demand $\bx$),
a price-taking supplier would like to employ the
scheduling policy $\pi$ and supply the bundle $\v{x}$.  We therefore have supply equal to demand.
 We show in the following corollary that the EDF scheduling policy $\pi^* = (\sigma^*, \phi)$ together with the
 marginal cost pricing scheme $\v{\zeta}$ constitute a market mechanism that induces a unique market equilibrium, which
maximizes the social welfare.

%
%\old{
%\begin{corollary}\label{coro:eff}
%Suppose that Assumptions \ref{a:meet}-\ref{A:IC} hold. For a given
%type distribution $\rho$, a mechanism $(\pi^*, \v{\zeta})$ maximizes
%the social welfare at the unique market equilibrium $(\bx,\v{\zeta}(\bx), \pi^*)$, where $\v{x}$ is the truth-telling aggregate demand bundle (cf. Eq. (\ref{eq:tt})) and
%$\pi^* = (\sigma^*,\phi)$ consists of an EDF inter-class scheduling policy $\sigma^*$ and an arbitrary intra-class policy $\phi \in \Phi(\sigma^*)$.
%\end{corollary}
%}

\begin{corollary}\label{coro:eff}
Suppose that Assumptions \ref{a:convex}-\ref{A:IC} hold. Given the
types of all consumers $\v{\theta}$ and  a market mechanism $(\pi^*,\v{\zeta})$,
there exists a unique market equilibrium
$(\bx^*,\v{\zeta}(\bx^*))$ that maximizes
the social welfare.
Here, $\bx^*$ is the truth-telling
aggregate demand bundle specified in  Equation \eqref{eq:demandB}.
\end{corollary}

\old{With the information on the aggregate demand, a regulatory authority
could induce the efficient equilibrium $(\bx^*,\v{\zeta}(\bx^*))$ by setting
the price bundle to be $\v{\zeta}(\bx^*)$. For deregulated markets,
rigorous treatment of the (dynamic) convergence to a competitive
equilibrium dates back to the early 20$^{\text{th}}$ century in the economics
literature \cite{H39,AH58,AH582}. For example, the instantaneous
adjustment process studied in \cite{AH58,AH582}, where the derivative of
price adjustment is set equal to the aggregate excess (net) demand
under the current price, is guaranteed to converge to the efficient competitive
equilibrium characterized in Corollary \ref{coro:eff}.}

%
%
%\old{With the information on consumers' aggregate demand, the supplier can simply announce the price bundle $\v{\zeta}(\bx)$, and maximize the social welfare
%at the competitive equilibrium $(\bx,\v{\zeta}(\bx))$.}

\section{Conclusion}\label{sec:con}

To explore the flexibility of deferrable electricity loads, we
propose a novel market for deadline-differentiated energy services that {offers
discounted (per-unit) electricity prices to consumers in exchange for their consent to defer
their electric power consumption, and}
provides a
guarantee on the aggregate quantity of energy to be delivered by a
consumer-specified deadline. We provide a full characterization of the joint scheduling and
pricing scheme that yields an efficient (competitive) market equilibrium between a (price-taking) supplier and a large consumer population. Somewhat
surprisingly, we show that this efficient mechanism is
incentive compatible in that every consumer would like to reveal her true
deadline to the supplier, regardless of the actions taken by other
consumers.

There are several interesting directions for future research.
First, it would be of practical interest to relax the assumption requiring a fixed price of firm supply  to allow for time-varying prices. 
Second, the market we have considered in our analysis is single shot, in the sense that it is cleared only once. As a natural extension, it would be of interest to analyze a 
dynamic analog of our formulation in which the market is
cleared on a dynamic basis over a finite or infinite horizon. Such a generalization of our model would also serve to facilitate the treatment of random consumer arrival times.

\begin{appendices}

%%%%%%%%%%%%%%%%%%%%%%%%%%%%%%%%%%%%%%%%%%%%%%%%%%%%%%%%%%%%%%%
\section{Intra-class Scheduling Policies} \label{app:intra}

Formally, for a consumer $i$ who purchases a bundle
$\v{a}=(a_1,\ldots,a_N)$, we let $ \v{\lambda}_{k,i}(\v{x}, \v{a})
\in \Rset_+^{N} $ denote (element-wise in $j$) the amount of energy
delivered (at period $k$) to consumer $i$ so as to satisfy her
demand $a_{j}$. We denote the intra-class scheduling policy by
\[ \phi = \left\{ (\v{\lambda}_{0,i}, \dots, \v{\lambda}_{N-1,i} )   \ | \   i \in [0,1]  \right\},\]
where $\v{\lambda}_{k,i} : \Rset_+^N \times \mathcal{A}
\rightarrow \Rset_+^N$ for all $i\in [0,1]$ and $k=0,\dots, N-1$.
Given an inter-class scheduling policy $\sigma \in \Sigma$, an intra-class
scheduling policy $\phi$ is feasible if and only if it satisfies the
following constraints.

\barablist \setlength{\itemsep}{.1in}
\item The intra-class scheduling policy should not deliver any
supply that is allocated to class $j$ to consumers outside this
class. That is, for each demand class $j=1,\ldots,N$, we have that
${\lambda}^j_{k,i}(\v{x}, \v{a})=0$ for every consumer $i$ such that
$a_{j}=0$.

\item At every period $k$, no energy is delivered to demand class
$j$ with $j \leq k$ by the feasibility of the inter-class policy $\sigma
\in \Sigma $. That is, $\lambda^{j}_{k,i}(\v{x}, \v{a})=0$ for every $i
\in [0,1]$, and every $1 \leq j \leq k \le N-1$.

\item The total supply allocated to demand class $j$ at time period $k$  must be fully utilized,
i.e.
\[
\int_{[0,1]} \lambda^{j}_{k,i}(\v{x}, \v{a}) \ \eta (di)  \ = \
\mu_k^j + \nu_k^j, \qquad \forall \ \ 0 \le k < j \le N-1,
\]
where we use the Lebesgue integral
%\footnote{Note that we have
%implicitly required here that $\lambda^{j}_{k,i}(\v{x}, \v{a})$ is
%Lebesgue integrable with respect to $i$, over the interval $[0,1]$.}
(with respect to Lebesgue measure $\eta$ defined on $[0,1]$), and
$\mu_k^j$ and $\nu_k^j$ denote the amount of intermittent and firm
supply allocated to demand class $j$ at period $k$, according to the
inter-class policy $\sigma$.

\item  Each consumer's individual delivery commitments must be met:
\[
\sum\nolimits_{t=1}^k a_{t} \ \leq  \ \omega_{k,i}^{\pi}(\bx,\v{a}) \ \le
\ \sum\nolimits_{t=1}^N a_{t}, \qquad  k=1,\ldots, N,
\]
where the total energy delivered to consumer $i$ by it's true deadline $k_i$ is
given by
\begin{equation}\label{eq:om}
 \omega_{k_i,i}^{\pi}(\v{x}, \v{a}) =  \sum_{t=0}^{k_i-1}  \sum_{j=1 }^N \lambda^j_{t,i}(\v{x}, \v{a}).
 \end{equation}
Notice that for any feasible inter-class scheduling policy $\sigma \in \Sigma$, it is always possible to ensure the satisfaction of the above constraint.
\end{list}
 We denote the \emph{feasible
intra-class policy space} by $\Phi(\sigma)$, which is
parameterized by a given inter-class policy $\sigma \in \Sigma$.

\section{Proof of Lemma \ref{lem:profit}} \label{sec:proof_lem}
We first establish the simplified form of $Q^*$ in Eq.  (\ref{eq:opt_val}). For notational simplicity, we denote the sequence of optimal inputs by $\v{u}_k^* = \v{\mu}_k^*(\v{z}_k, s_k)$ and $\v{v}_k^* = \v{\nu}_k^*(\v{z}_k, s_k)$ for $k=0,\dots, N-1$. Under the optimal scheduling policy, firm supply is deployed only as a last resort to ensure task satisfaction. It follows that
\[ Q^*(\v{x}) \ = \ c_0 \cdot \sum_{k=1}^N \mathbb{E} \left\{ v_{k-1}^{k,*} \right\}. \]
To establish the desired result, it suffices to show that $v_{k-1}^{k,*} = - \min\{0,\xi_k\}$. First, define the quantity
\[
 \delta_k \ = \ \sum_{j=0}^{k-1} s_j   \ - \ \sum_{\ell=1}^k \sum_{j=0}^{\ell-1} u_j^{\ell,*}, \qquad \forall \ \   k=1,\dots,N,
 \]
which denotes the maximum amount of intermittent supply available to demand class $k+1$  across the first  $k$ time periods under a sequence of EDF allocations $\v{u}^*_0, \dots, \v{u}^*_{k-1}$. Clearly, we have that $\delta_k \geq 0$ for all $k$, given feasibility of the allocations $\v{u}^*_0, \dots, \v{u}^*_{k-1}$ under  the intermittent supply availability constraints. One can readily show via an inductive argument that
\[ \xi_k \ = \ \delta_k \ - \ v_{k-1}^{k,*},   \qquad \forall \ \   k=1,\dots,N. \]
Using this characterization of the residual process $\v{\xi}$, we have that
\begin{align*}
&\delta_k > 0 \ \ \Longrightarrow \  \ v_{k-1}^{k,*} = 0  \ \ \Longrightarrow \  \  \min\{0,\xi_k\} = 0,\\
&\delta_k = 0  \ \ \Longrightarrow \  \  v_{k-1}^{k,*}  \geq 0 \ \ \Longrightarrow \  \  \min\{0,\xi_k\} = -v_{k-1}^{k,*},
\end{align*}
which yields the desired result that  $\min\{0,\xi_k\} \ = \  - \ v_{k-1}^{k,*}$ and establishes the form of $Q^*$ in Eq.  (\ref{eq:opt_val}).

We now establish \emph{convexity} of the expected recourse cost  $Q^*(\v{x})$ directly.
Let $\v{x_1}, \v{x_2} \in \Rset_+$ and denote the corresponding optimal scheduling policies by $\sigma_1^* \in \Sigma(\v{x}_1)$ and $\sigma_2^* \in \Sigma(\v{x}_2)$.  Define the convex combination of demand bundles
$ \v{x}_{\lambda} = \lambda \v{x}_1 + (1-\lambda) \v{x}_2$,
where $\lambda \in [0,1]$. It follows that
\begin{align*}
\lambda Q^*(\v{x}_1) + (1-\lambda) Q^*(\v{x}_2)   \ = \ \mathbb{E} \sum_{k=0}^{N-1}    \v{c}^{\top}( \lambda \v{v}_k^{\sigma_1^*}  + (1-\lambda) \v{v}_k^{\sigma_2^*}).
\end{align*}
And, it is not difficult to show that the convex combination of the constituent policies $\lambda \sigma_1^*  + (1-\lambda) \sigma_2^*$ is admissible for the convex combination of demand bundles, i.e.,
$ \lambda \sigma_1^*  + (1-\lambda) \sigma_2^* \in \Sigma(\v{x}_{\lambda})$. Convexity of $Q^*$ follows.
\emph{Differentiability} of  $Q^*(\v{x})$ over $(0,\infty)^N$ follows immediately from the proof of Theorem \ref{thm:opt_sup}, in which we show that
$$ \frac{\partial Q^*(\v{x})}{\partial x_k}  =  \zeta_k(\v{x}), \quad k=1,\dots, N,$$
 where the function $\zeta_k(\v{x})$, defined in \eqref{eq:eP}, is  bounded and continuous over $(0,\infty)^N$ for each $k=1,\dots, N$.  $\hfill \blacksquare$

%%%%%%%%%%%%%%%%%%%%%%%%%%%%%%%%%%%%%%%%%%%%%%%%%%%%%%%%%%%%%%%
\section{Proof of Theorem \ref{thm:opt_sup}} \label{ap:opt_sup}
To simplify exposition, we employ the shorthand notation $(\cdot)^- = \min\{0, \cdot\}$.
Fix a price bundle $\v{p} \in \Rset_+^N$. We have previously shown in Theorem \ref{thm:opt_pol_app} that the EDF inter-class scheduling policy $\sigma^* \in \Sigma(\v{x})$ is optimal for any demand bundle $\v{x} \in \Rset_+^N$.
Hence, it suffices to show that $\v{p} = \v{\zeta}(\v{x})$ satisfies the first order condition for optimality (\ref{eq:foc}).
Taking the gradient of the supplier's expected profit with respect to $\v{x}$ yields \
$\nabla_{\v{x}} \;  J(\v{x}, \v{p}, \sigma^*) \ = \ \v{p}  \ - \ \nabla_{\v{x}} \; Q^*(\v{x}).$
It remains to show that
\[ \frac{\partial Q^*(\v{x})}{\partial x_k}  \ = \   \zeta_k(\v{x})\]
for $ k=1,\dots, N.$ Working the with the simplified form of $Q^*$ established in Lemma \ref{lem:profit},
it follows readily  that
\begin{align} \label{eq:id1}
\frac{\partial Q^*(\v{x})}{\partial x_k} \ = \- c_0 \cdot \sum_{\ell=k}^N \frac{\partial}{\partial x_k} \; \mathbb{E} \left\{ \xi_{\ell}(\v{x}, \v{s})^-  \right\},
\end{align}
where we've truncated the summation from below at $\ell=k$,\ as \ $\xi_{\ell}(\v{x}, \v{s})$ is wholly independent of $x_k$ for all $\ell  < k$.  We therefore restrict our attention to $\ell \geq k$  for the remainder of the proof.
The next step of the proof relies on the ability to interchange the order of differentiation and expectation in (\ref{eq:id1}). It is obvious from construction that $\xi_{\ell}(\v{x}, \v{s})^-$ is both a continuous function of $(\v{x}, \v{s})$ and piecewise affine in $x_k$ (with a finite number of linear segments) for each $\v{s}$.
It follows that $\xi_{\ell}(\v{x}, \v{s})^-$  is differentiable almost everywhere in $x_k \in \Rset_+$ and satisfies
\[ \left| \frac{\partial}{\partial x_k}  \xi_{\ell}(\v{x}, \v{s})^-   \right| \  \leq \ 1   \]
almost everywhere.  Then, for each $x_k \in \Rset_+$, we have that  $\partial \xi_{\ell}(\v{x}, \v{s})^-/ \partial x_k$ is integrable in $x_k$ and by the dominated convergence theorem
\begin{align*}
\frac{\partial}{\partial x_k} \; \mathbb{E} \left\{ \xi_{\ell}(\v{x}, \v{s})^-  \right\} & \ = \  \mathbb{E} \left\{ \frac{\partial}{\partial x_k}  \xi_{\ell}(\v{x}, \v{s})^-  \right\}.
\end{align*}
Finally, it is not difficult to see that
\begin{align*}
\frac{\partial}{\partial x_k}  \xi_{\ell}(\v{x}, \v{s})^-  \ = \
\left\{
\begin{array}{ll}
\displaystyle \v{1}_{ \{  \xi_k \leq 0 \} }, & \quad \ell = k, \\[8pt]
\displaystyle \v{1}_{ \{  \xi_{\ell} \leq 0 \} } \cdot \prod_{t=k}^{\ell-1} \  \v{1}_{ \{  \xi_{t} > 0 \} } , & \quad \ell > k.
\end{array} \right.
\end{align*}
And taking expectation, we have the desired result. $\hfill \blacksquare$

\section{Proof of Theorem \ref{t:IC}}\label{sec:A}

 \old{We first note
that for any $\v{x}$, the price bundle $\v{\zeta}(\v{x})$ (resulting
from the price schedule in \eqref{eq:eP}) is nonincreasing in
 the deadline $k$, and that $\zeta_k(\v{x}) \le c_0 $ for every
 $k$.}

 Let $\bx$ be the aggregate
demand of the other consumers (excluding $i$). Suppose that the supplier uses the EDF (inter-class scheduling) policy
$\sigma^*$, and an arbitrary intra-class scheduling policy $\phi \in \Phi(\sigma^*)$. We let $\pi^*=(\sigma^*,\phi)$.
We consider a consumer $i$ of some type $\theta_i=(k,R,q)$ such that
$c_0 \le R$. We will show that the consumer, who faces an arbitrary aggregate demand bundleh, $\bx$,
requested by other consumers, would like to take the truth-telling action specified in Definition \ref{def:truth}.
For the rest of this proof, we will use
$\v{p}=\{p_k\}_{k=1}^N$ to denote the price bundle induced by the pricing
scheme $\v{\zeta}$ (cf. its definition in Eq. \eqref{eq:eP}), at the
aggregate demand $\bx$.

If consumer $i$ is
truth-telling, she will request a quantity $q$ before deadline $k$,
and receive an expected payoff of
\begin{equation}\label{eq:pt}
\begin{array}{l}
 \displaystyle \;\;\;\; V^{\pi^*}_{i}(\theta_i,\varphi^*_i, \bx,\v{\zeta})\\[5pt]
 \displaystyle = U(q)-qp_k \\[5pt]
\displaystyle = qR-q c_0 \Big[
\mathbb{P}(\xi_k \le 0)  \\[5pt]
 \displaystyle  \quad \qquad\qquad
\displaystyle 
 + \sum\nolimits_{t=k+1}^N
\mathbb{P}(\xi_k
> 0,\ldots,\xi_{t-1}>0, \xi_{t}\le 0)\Big].
\end{array}
\end{equation}

Since the optimal price schedule is nonincreasing in deadline, and
demand is guaranteed to be met before the requested deadline, the
consumer has no incentive to request a positive amount of
electricity at some period $t$ that is earlier than $k$. We can
therefore assume that consumer $i$ takes an action
$\v{a}'=\varphi'_i(\theta_i)$ such that
$$
a'_t =0,\qquad t=0,\ldots,k-1.
$$

Since the consumer cannot increase its expected payoff \blue{(compared to being truth-telling)} by reporting
$a'_k \ge q$, we focus on the case where $a'_k < q$. We first write
the consumer's expected payoff (achieved by the action $\v{a}'$) as
$$
V^{\pi^*}_{i}(\theta_i,\varphi'_i, \bx,\v{\zeta}) =
 \E \left\{U_\theta\left( \omega^{\pi^*}_{k,i}(\bx,\v{a}')\right)\right\} - \sum\limits_{t=k}^N p_t
a'_t.
$$
Showing that $V^{\pi^*}_{i}({\theta}_i,\varphi'_i,
\bx,\v{\zeta})$ is no more than the expected payoff in
\eqref{eq:pt} is equivalent to
\begin{equation}\label{eq:sh}
\begin{array}{l}
\displaystyle \;\;\;\; R q- (q-a'_k)c_0 \Big[
\mathbb{P}(\xi_k \le 0) \\[5pt]
\displaystyle \qquad \qquad \;\;\;\;
 +\sum\nolimits_{t=k+1}^N \mathbb{P}(\xi_k
> 0,\ldots,\xi_{t-1}>0, \xi_{t}\le
0)\Big] \\[8pt]
\displaystyle \ge \E \left\{U_\theta\blue{\left(\omega^{\pi^*}_{k,i}(\bx,\v{a}')
\right)}\right\}- \sum\nolimits_{t=k+1}^N p_t a'_t.
\end{array}
\end{equation}
 We will derive an
upper bound on the right hand side of \eqref{eq:sh}, and show that
this upper bound cannot exceed the left hand side of \eqref{eq:sh}.
For notational convenience, we define
$$
\eta_t = \{\xi_k>0, \ldots, \xi_{t-1}>0,\xi_t \le0 \}
$$
for all $t=k+1,\ldots,N-1$ and
$$
\blue{\eta_k=\{\xi_k \le 0\}}, \qquad \eta_{N} = \{\xi_k>0, \ldots, \xi_{N-1}>0 \}.
$$
Note that these $N-k+1$ events are mutually disjoint. Further, it is
straightforward to see that $\xi_k \le 0$ implies that
$\omega^\pi_{k,i}(\bx,\v{a}')=a'_k$. While, on the other hand, $\xi_k
> 0$ implies that one of the (mutually disjoint) events
$\{\eta_t\}_{t=k+1}^{N}$ must occur, i.e.,
\begin{equation}\label{eq:11}
1=\mathbb{P} \left( \bigcup\nolimits_{t=k}^{N}
\eta_t \right)= \sum_{t=k}^{N}
\mathbb{P}(\eta_t ).
\end{equation}
\blue{Under the event $\eta_t$,
the amount of energy delivered by the EDF scheduling policy before her true deadline $k$,
$\omega^{\pi^*}_{k,i}(\bx,\v{a}')$, cannot exceed $\sum\nolimits_{\tau=k}^t a'_{\tau}$, for $t=k,\ldots,N$.
It follows from Assumption \ref{a:convex} that}
$$
\blue{ Rq- U_\theta\left( x \right) \ge R \left(q - x \right)^+, \qquad \forall\; x \ge 0 },
$$
where $(\cdot)^+=\max\{\cdot,0\}$. We then have
\begin{equation}\label{eq:sh1}
\begin{array}{l}
\displaystyle\;\;\;\; Rq-\E \left\{U_\theta\left(
\omega^{\pi^*}_{k,i}(\bx,\v{a}')\right)\right\} \\ [5pt]
\displaystyle \ge  R\sum\nolimits_{t=k}^{N} \mathbb{P}(\eta_t) \cdot \left(q - \sum\nolimits_{\tau=k}^t a'_{\tau}  \right)^+.
\end{array}
\end{equation}

 We now argue that the right hand side of \eqref{eq:sh1} is
minimized at some vector $\tilde {\v{a}}$ such that $\sum_{t=k}^N
\tilde a_t \le q$. To see this, suppose that $\sum_{t=k}^N a'_t >
q$. Let $T$ be the smallest $t$ such that $\sum_{m=k}^t  a'_m > q$.
We define an alternative vector $\tilde {\v{a}}$
\begin{equation}\label{eq:aa}
\tilde a_t=a'_t,\qquad t \le T-1,\qquad \tilde a_T= q-
\sum_{m=k}^{T-1}  a'_m,
\end{equation}
and  $\tilde a_t=0$, for every $t>T$.
 We have $\sum_{t=k}^N \tilde a_t = q$, and
$$
\begin{array}{l}
\displaystyle\;\;\;\; 
\sum\nolimits_{t=k}^{N} \mathbb{P}(\eta_t) \cdot \left(q - \sum\nolimits_{\tau=k}^t a'_{\tau}  \right)^+  \\[5pt]
\displaystyle = \sum\nolimits_{t=k}^{N}
\mathbb{P}(\eta_t) \cdot \left(q - \sum\nolimits_{\tau=k}^t \tilde a_{\tau}  \right)^+.
\end{array}
$$
Note that if $\sum_{t=k}^N a'_t \le q$, then $\tilde
{\v{a}}=\v{a'}$.
  To validate
\eqref{eq:sh}, it suffices to show that for any vector $\tilde
{\v{a}}$ defined in \eqref{eq:aa},
\begin{equation}\label{eq:sh2}
\begin{array}{l}
\displaystyle \;\;\;\; \blue{R\sum\nolimits_{t=k}^{N}
\mathbb{P}(\eta_t) \cdot \left(q - \sum\nolimits_{\tau=k}^t \tilde a_{\tau}  \right) }\\[13pt]
\displaystyle \ge - \sum\nolimits_{t=k+1}^N p_t \tilde a_t + (q- \tilde a_k)c_0 \Big[
\mathbb{P}(\xi_k \le 0) \\[8pt]
\displaystyle \qquad \;\;\;\;\;\qquad +\sum\nolimits_{t=k+1}^N \mathbb{P}(\xi_k
> 0,\ldots,\xi_{t-1}>0, \xi_{t}\le
0)\Big],
\end{array}
\end{equation}
\blue{where the left hand side is a lower bound on the loss of expected utility due to
the non-truthful action $\tilde {\v{a}}$ (cf. the inequality in \eqref{eq:sh1}),
and the right hand side is the difference between the payment under the truth-telling action and the action
$\tilde {\v{a}}$. We will prove the following inequality that is equivalent to \eqref{eq:sh2}}
\begin{equation}\label{eq:sh20}
\begin{array}{l}
\displaystyle \;\;\;\; \blue{R\sum\nolimits_{t=k}^{N}
\mathbb{P}(\eta_t) \cdot \left(q - \sum\nolimits_{\tau=k}^t \tilde a_{\tau}  \right) }\\[13pt]
\displaystyle =R(q-\tilde a_k) - R \sum\nolimits_{t=k+1}^{N}
\left( \mathbb{P}(\eta_t) \cdot \sum\nolimits_{\tau=k+1}^t \tilde a_{\tau}  \right)\\[13pt]
\displaystyle =R(q-\tilde a_k) - R \sum\nolimits_{t=k+1}^{N}
\left(\tilde a_{t} \cdot \sum\nolimits_{\tau=m}^N  \mathbb{P}(\eta_m) \right)\\[13pt]
\displaystyle \ge (q- \tilde a_k)c_0 \Big[
\mathbb{P}(\xi_k \le 0)\\[7pt]
\displaystyle  +\sum\limits_{t=k+1}^N \mathbb{P}(\xi_k
> 0,\ldots,\xi_{t-1}>0, \xi_{t}\le
0)\Big]- \sum\nolimits_{t=k+1}^N p_t \tilde a_t,
\end{array}
\end{equation}
\blue{where the first equality follows from \eqref{eq:11}, and the second equality is obtained by rearranging terms.}

It follows from the characterization of supplier marginal cost
in \eqref{eq:eP}
that for $t=k+1,\ldots,N,$
\begin{equation}\label{eq:u2}
\begin{array}{l}
\displaystyle\;\;\;\; \tilde a_{t} p_{t} \ = \ \tilde a_{t} c_0 \sum_{m=t}^{N}
\mathbb{P}(\xi_t
> 0,\ldots,\xi_{m-1}>0, \xi_{m}\le 0)  \\[6pt]
\displaystyle \ge \ \tilde a_{t} c_0 \sum_{m=t}^{N}
\mathbb{P}(\xi_k
> 0,\ldots,\xi_{m-1}>0, \xi_{m}\le 0) .
\end{array}
\end{equation}
We then have, for $t=k+1,\ldots,N$,
\begin{equation}\label{eq:sh22}
\begin{array}{l}
\displaystyle \;\;\;\; R \tilde a_t \sum\nolimits_{m=t}^{N}
\mathbb{P}(\eta_m) - \tilde a_t p_t \\[8pt]
\displaystyle = R \tilde a_t \left( 1- \mathbb{P} (\xi_k \le 0) - \sum\nolimits_{m=k+1}^{t-1} \mathbb{P} (\eta_m) \right) - \tilde a_t p_t \\[8pt]
\displaystyle \le R \tilde a_t \left( 1- \mathbb{P} (\xi_k \le 0) -
\sum\nolimits_{m=k+1}^{t-1} \mathbb{P} (\eta_m) \right) \\[6pt]
\displaystyle \qquad \qquad
- \tilde a_{t} c_0 \sum\nolimits_{m=t}^{N}
\mathbb{P}(\xi_k> 0,\ldots,\xi_{m-1}>0, \xi_{m}\le 0)\\[8pt]
\displaystyle  \le R \tilde a_t \left( 1- \mathbb{P} (\xi_k \le 0)
\right) \\[6pt]
\displaystyle \qquad  \;- \tilde a_{t} c_0 \sum\nolimits_{m=k+1}^{N}
\mathbb{P}(\xi_k
> 0,\ldots,\xi_{m-1}>0, \xi_{m}\le 0)\\[8pt]
\displaystyle \le R \tilde a_t - \tilde a_{t} c_0  \Big(\mathbb{P} (\xi_k \le 0) \\[6pt]
\displaystyle \qquad \qquad+ \sum\nolimits_{m=k+1}^{N}
\mathbb{P}(\xi_k
> 0,\ldots,\xi_{m-1}>0, \xi_{m}\le 0)\Big).
\end{array}
\end{equation}
Here, the first inequality follows from \eqref{eq:u2}; the second
inequality is true, because $R \ge c_0$ and $\eta_m = \{\xi_k
> 0,\ldots,\xi_{m-1}>0, \xi_{m}\le 0\}$, for $m=k+1,\ldots,t-1$;
the last inequality follows from $R \ge c_0$. \blue{For any vector $\tilde
{\v{a}}$ with $\sum\nolimits_{t=k}^N \tilde a_t \le q$},
from \eqref{eq:sh22}
we have
\begin{equation}\label{eq:sh3}
\begin{array}{l}
\displaystyle \;\;\;\;  \sum\nolimits_{t=k+1}^N \left(-p_t \tilde
a_t+
R \tilde a_t \sum\nolimits_{m=t}^{N} \mathbb{P}(\eta_m) \right)\\[13pt]
\displaystyle \le \sum\nolimits_{t=k+1}^N \tilde a_t \Big[ R-c_0
\Big(\mathbb{P} (\xi_k \le 0)   \\[8pt]
\displaystyle \qquad \qquad
+\sum\nolimits_{m=k+1}^{N}
\mathbb{P}(\xi_k
> 0,\ldots,\xi_{m-1}>0, \xi_{m}\le 0)\Big) \Big].
\end{array}
\end{equation}
Since $R \ge c_0$, the right hand side of \eqref{eq:sh3} is
nondecreasing in $\sum\nolimits_{t=k+1}^N \tilde a_t$. The desired
result in \eqref{eq:sh20} follows from the fact
$\sum\nolimits_{t=k+1}^N \tilde a_t \le q-\tilde a_k$.
 Since the preceding analysis holds for any action $\v{a}'$ and any aggregate demand bundle $\bx$, we conclude
 that it is a dominant strategy for consumer $i$ to be
 truth-telling, i.e., the
 pricing scheme \eqref{eq:eP} is incentive compatible, in the sense of Definition \ref{Def:ICM}.  $\hfill \blacksquare$

\section{Proof of Corollary \ref{coro:eff}}\label{sec:D}

We first show that $(\bx^*,\v{\zeta}(\bx^*))$ is a market equilibrium (in the sense of Definition \ref{def:CE}),
and then argue that it is the unique market equilibrium that maximizes the social welfare.
Under the mechanism $(\sigma^*,\phi,\v{\zeta})$, it follows from Assumption \ref{A:IC} and Theorem \ref{t:IC}
that it is a dominant strategy for every consumer to be
 truth-telling, and therefore the aggregate demand is  $\bx^*$.
 Given the price bundle $\v{\zeta}(\bx^*)$,
 the aggregate demand bundle $\v{x}^*$ together
 with the EDF scheduling policy maximizes the supplier's
expected profit (cf. Theorem \ref{thm:opt_sup}). Hence,
 the pair, $(\v{x}^*, \v{\zeta}(\v{x}^*))$, constitutes a market
 equilibrium.

 Let $(\bx, \v{p})$ be a market equilibrium.
 The second condition of Definition \ref{def:CE}
 requires that the given the price bundle $\v{p}$,
 the quantity $\bx$ maximizes the (price-taking) supplier's expected profit.
 We note that this implies that $\v{p}$ is the supplier's marginal cost to supply
 $\bx$. Since the supplier marginal cost never exceeds $c_0$, at any price bundle that represents supplier marginal cost,
 a truth-telling consumer would request her maximum demand by her true deadline, and therefore
  the aggregate demand of a truth-telling consumer population  
  is always $\bx^*$. The uniqueness of the truth-telling aggregate demand
implies the uniqueness of a market equilibrium.

It is straightforward to see that social welfare is maximized at this market equilibrium, because
under Assumption \ref{A:IC}, it is socially optimal to fully serve the aggregate demand $\v{x}^*$, and further, the EDF scheduling policy $\sigma^*$
minimizes the expected cost of servicing $\bx^*$. $\hfill \blacksquare$

\end{appendices}

%%%%%%%%%%%%%%%%%%%%%%%%%%%%%%%%%%%%%%%%%%%%%%%%%%%%%%%%

\begin{IEEEbiography}[{\includegraphics[width=1in,height=1.25in,clip,keepaspectratio]{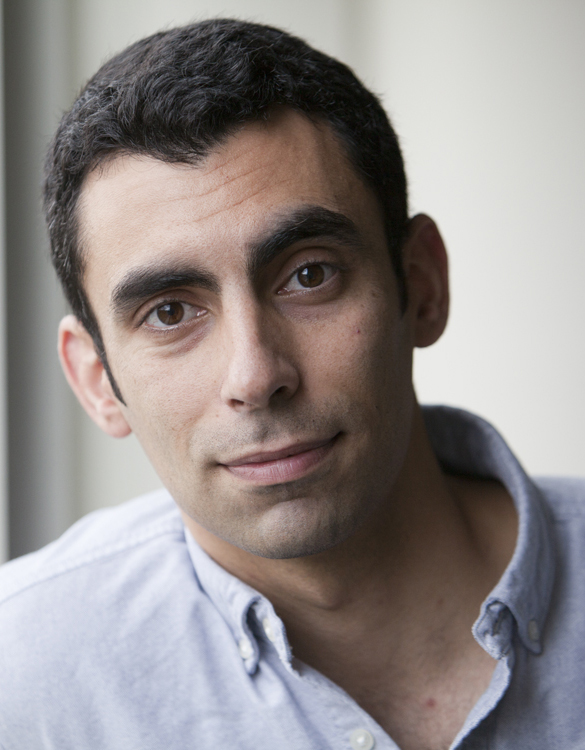}}]{Eilyan Bitar}
currently serves as an Assistant Professor and the David D. Croll Sesquicentennial Faculty Fellow in the School of Electrical and Computer Engineering at Cornell University, Ithaca, NY, USA.  Prior to joining Cornell in the Fall of 2012, he was engaged as a Postdoctoral Fellow in the department of Computing and Mathematical Sciences at the California Institute of Technology and at the University of California, Berkeley in Electrical Engineering and Computer Science, during the 2011-2012 academic year. His current research examines the operation and economics of modern power systems, with an emphasis on the design of markets and optimization methods to manage uncertainty in renewable energy resources. He received the B.S. and Ph.D. degrees in Mechanical Engineering from the University of California at Berkeley in 2006 and 2011, respectively.

Dr. Bitar is a recipient of the NSF Faculty Early Career Development Award (CAREER), the John and Janet McMurtry Fellowship, the John G. Maurer Fellowship, and the Robert F. Steidel Jr. Fellowship.
\end{IEEEbiography}

\begin{IEEEbiography}[{\includegraphics[width=1in,height=1.25in,clip,keepaspectratio]{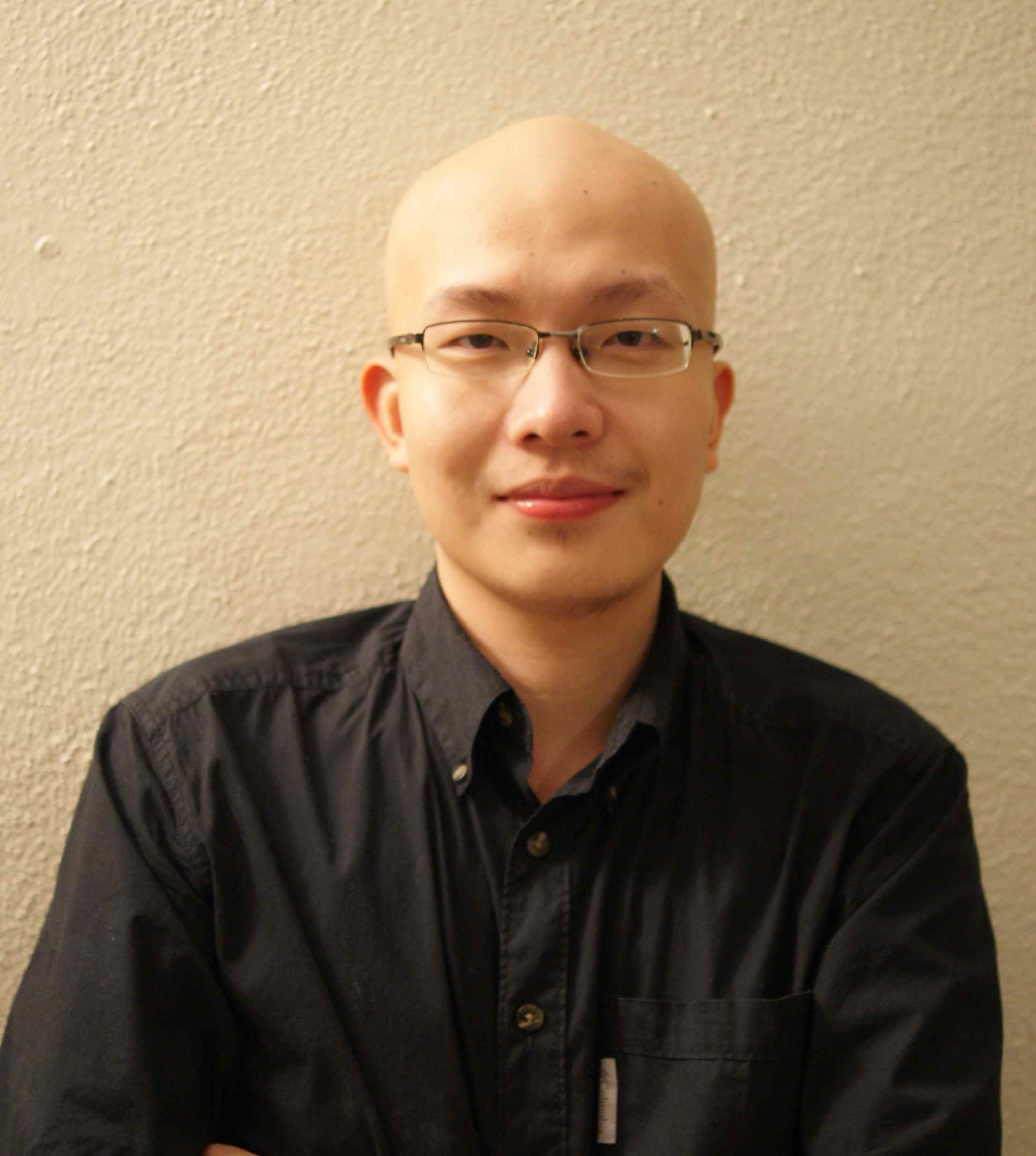}}]{Yunjian Xu} (S'06-M'10) received the B.S. and M.S. degrees in
electrical engineering from Tsinghua
University, Beijing, China, in 2006 and 2008, respectively, and the
Ph.D. degree from the Massachusetts Institute of
Technology (MIT), Cambridge, MA, USA, in 2012.

Dr. Xu was a CMI (Center for the Mathematics of Information)
postdoctoral fellow at the California Institute of Technology,
Pasadena, CA, USA, in 2012-2013.
He joined the Singapore University of Technology and Design,
Singapore, as an assistant professor in 2013. His research interests
focus on energy systems and markets, with emphasis on power system
optimization, wholesale electricity market
design, and the aggregation of distributed energy resources in power
distribution systems.

Dr. Xu was a recipient of the MIT-Shell Energy Fellowship.
\end{IEEEbiography}

%%%%%%%%%%%%%%%%%%%%%%%%%%%%%%%%%%%%%%%%%%%%%%%%%%%%%%%%%%%%%%%%%%%%%%%%%

\end{document}